\documentclass[a4paper, 10pt]{article}
\usepackage{cite}
\usepackage{amsfonts}
\usepackage{amscd}
\usepackage{amsmath}
\usepackage{amsthm}
\usepackage[all]{xy}
\usepackage[dvips]{graphicx}
\usepackage{verbatim}
\usepackage{enumerate}
\usepackage{amsmath}
\usepackage{amssymb}
\usepackage{url}
\usepackage{xcolor}
\usepackage{mathrsfs}
\usepackage{authblk}
\usepackage{color, xcolor}
\usepackage{psfrag}
\usepackage{comment}
\usepackage{tikz}
\usepackage{hyperref}
\usepackage{indentfirst}
\usepackage{comment}

\usepackage[margin=1.2in]{geometry}

\newcommand{\RR}{\mathbb{R}}

\newcommand{\CC}{\mathbb{C}}
\newcommand{\NN}{\mathbb{N}}
\newcommand{\CP}{\mathcal{P}}

\newcommand{\ZZ}{\mathbb{Z}}
\newcommand{\TT}{\mathbb{T}}

\newcommand{\CB}{\mathcal{B}}

\newcommand{\CDD}{\mathcal{D}}
\newcommand{\CE}{\mathcal{E}}

\newcommand{\CU}{\mathcal{U}}
\newcommand{\CV}{\mathcal{V}}
\newcommand{\CL}{\mathcal{L}}
\newcommand{\CJ}{\mathcal{J}}

\newcommand{\CN}{\mathcal{N}}

\newcommand{\CR}{\mathcal{R}}

\newcommand{\ga}{\gamma}
\newcommand{\al}{\alpha}

\newcommand{\p}{\partial}

\newcommand{\wrho}{\widetilde{\rho}}

\newcommand{\hv}{\widehat{v}}

\newcommand{\lan}{\langle}
\newcommand{\ran}{\rangle}

\newcommand{\ov}{\overline{v}}

\newcommand{\beq}{\begin{equation}}
\newcommand{\eeq}{\end{equation}}
\newcommand{\beqn}{\begin{equation*}}
\newcommand{\eeqn}{\end{equation*}}
\newcommand{\beqy}{\begin{eqnarray}}
\newcommand{\eeqy}{\end{eqnarray}}
\newcommand{\beqyn}{\begin{eqnarray*}}
\newcommand{\eeqyn}{\end{eqnarray*}}

\newtheorem{theorem}{Theorem}[section]
\newtheorem{prop}[theorem]{Proposition}
\newtheorem{lemma}[theorem]{Lemma}

\newtheorem{cor}[theorem]{Corollary}
\newtheorem{remark}[theorem]{Remark}
\newtheorem{definition}[theorem]{Definition}

\newtheorem{example}[theorem]{Example}

\newtheorem{problem}[theorem]{Problem}

\newtheorem{propertyx}{Property}

\numberwithin{equation}{section}

\newcommand{\ds}{\displaystyle}

\newcommand{\Par}{\mathrm{Par}}
\newcommand{\bi}{\mathbf{i}}
\newcommand{\bj}{\mathbf{j}}

\newcommand{\EEnd}{\mathrm{EEnd}}

\newcommand{\bv}{\mathbf{v}}
\newcommand{\dive}{\mathrm{div}}

\newcommand{\ww}{\widetilde{w}}

\newcommand{\rct}{\mathrm{ct}}

\title{Explicit Response Formulae and Request Problem \\
 for Smooth Expanding Endomorphisms}

\author[1]{Jianyu Chen\thanks{jychen@suda.edu.cn}}
\author[1]{Na Zhang\thanks{ 3199885810@qq.com}}
\author[2]{Yufeng Zhang\thanks{zyfp@mail.ustc.edu.cn}}
\affil[1]{\footnotesize School of Mathematical Sciences
\& Center for Dynamical Systems and Differential Equations,
\hspace{3cm}
Soochow University, Suzhou, Jiangsu, P.R.China}
\affil[2]{\footnotesize Mathematics \& Statistical Sciences,
University of Alberta, Edmonton, Canada}

\date{}

\begin{document}

\maketitle

\begin{abstract}

In this paper, 
we investigate the higher order differentiability of physical measures,
in the category of smooth expanding endomorphisms of compact manifolds,
and we establish the explicit response formulae of any order.
The proof is based on the elegant Keller-Liverani-Gou\"ezel perturbation theory of transfer operators. 
Furthermore, we provide the solutions to the high order request problem,
which is regarded as the inverse problem to response formulae,
that is, we make perturbations to realize prescribed responses up to any finite order, and we establish the existence of minimizers for a relaxed optimization problem associated with suitable cost functionals.
\\


\end{abstract}

\tableofcontents

\medskip

\section{Introduction}

The starting point of this paper is to investigate 
the stability theory of hyperbolic dynamical systems
from the viewpoint of ``physical measures'',
which include the volume measure of conservative systems
and 
the Sinai-Ruelle-Bowen (SRB) measure (an analogy of volume) of 
dissipative systems~\cite{MR0399421, MR380889, MR0442989, MR415683}.
In particular, the physical measures for uniformly expanding systems are 
in fact absolutely continuous invariant probability measures (ACIPs) with respect to the volume. 
Perturbing the dynamics may change its physical measure and hence
its long-term behavior.
Understanding these changes
leads naturally to the {\it response problem}: how do physical measures
depend on the perturbation, and can their infinitesimal changes
be expressed explicitly?
Conversely, one may seek perturbations that produce prescribed
changes in physical measures.
This inverse problem is called the {\it request problem}.
These two problems play an important role in applications, 
including the assessment of climate sensitivity~\cite{MR4166837}
and the control of long-term statistics of dynamical
systems~\cite{MR3670005,MR3787356}.
 
\medskip

To formulate the stability of physical measures,
we consider a continuous   family of dynamical systems $\{f_t\}$ on a phase space $X$, 
where the parameter $t\in (-t_0, t_0)$ for some $t_0>0$, and 
we assume that 
there is a unique physical measure $\mu_t$ for each system $f_t$. 
It is natural to ask the following question:
is the mapping $t\mapsto \mu_t$ weakly continuous? That is, 
given a smooth  function $\psi$ on $X$, 
is the mapping $\ds t\mapsto \int \psi \mu_t$ continuous?
It is worth pointing out that such weak continuity   
might fail for non-smooth hyperbolic systems, for instance, Keller~\cite{MR685377} and later Li et al~\cite{MR3009107} showed the instability of ACIPs for a continuous family of Keller's $W$-maps. 
 
\medskip 
 
We restrict ourselves on smooth parameter families of smooth dynamical systems. 
The pioneering work of Ruelle~\cite{MR1463827} achieved not only the continuity but also 
the differentiability of SRB measures $\{\mu_t\}$ on mixing Axiom A attractors of $\{f_t\}$.
Moreover, the following linear response formula was established:
\beqn
\p_t\int\psi d\mu_t 
=\sum_{n=0}^\infty \int \lan \nabla ( \psi\circ f_t^n ), V_t \ran \ d\mu_t,
\eeqn
where the vector field $V_t$ is the infinitesimal variation of $\{f_t\}$, 
i.e., $\ds \left.\p_s f_s\right|_{s=t}= V_t\circ f_t$. 
A similar result was later obtained by Ruelle for the hyperbolic flows~\cite{MR2408395}.
In the particular case when
$\{f_t\}$ are smooth expanding circle maps, 
Baladi~\cite{MR3729040} gave an alternative expression of the linear response formula
in terms of ACIP densities $\ds \rho_t(x)=d\mu_t(x)/dx$, that is, 
\beqn
\p_t \rho_t  = -\left(I-\CL_t\right)^{-1} \p_x\left(\rho_tV_t\right),
\eeqn
where $\CL_t$ is the transfer operator of $f_t$ with respect to Lebesgue measure.
These classical results motivate us to further explore the following two directions:
\begin{enumerate}
\item[(i)] 
Does the linear response formula hold 
for dynamical systems beyond uniform hyperbolicity 
(including nonuniformly partially hyperbolic systems, possibly with singularities)?
\item[(ii)] 
With sufficient smoothness, do the physical measures $\{\mu_t\}$ exhibit higher order  differentiability, i.e., 
can $\ds t\mapsto  \mu_t$ be $C^k$ weakly differentiable for $k\ge 2$? 
Furthermore, can we obtain explicit response formulae of higher order?
\end{enumerate}

\medskip

Developments in direction (i) were mostly   for 
concrete systems under fairly strong conditions, 
for instance, Dolgopyat~\cite{MR2031432} established  
that the linear response of SRB measures holds for rapidly mixing abelian Anosov actions;
however, Zhang~\cite{MR3772032} constructed a   
partially hyperbolic diffeomorphism on the 3-torus at which the differentiability of SRB measures fails;
Baladi and Smania obtained the linear response of ACIPs for
piecewise expanding unimodal maps within tangency of conjugacy classes~\cite{MR2399821},
and also for the S-unimodal maps satisfying the topologically slow recurrent condition~\cite{MR3075107};
the linear response for the  Pomeau-Manneville intermittent maps (with respect to exponent parameter)
was shown by 
Baladi and Todd~\cite{MR3551256} and by Korepanov~\cite{MR3502226} independently. 
Surprisingly, the linear response has also recently been established
for certain systems with singularities.
For instance, Bahsoun and Galatolo~\cite{MR4752563} proved linear
response for cusp maps, and Canestrari~\cite{MR5068786}
established it for discontinuous perturbations of smooth endomorphisms.

\medskip

Regarding direction (ii), 
Ruelle~\cite{MR1492947} gave a formal rather than rigorous derivation for the 
higher order response formulae. 
Recently, the quadratic response formula was obtained by Galatolo and  Sedro~\cite{MR4059336}, and also by
Galatolo and Marsan~\cite{MR4305587}, in both random and deterministic dynamical systems;
Crimmins and Nakano~\cite{MR4713046} obtained
quenched linear and higher order response for some partially hyperbolic systems. 
Higher order differentiability is also treated by the general perturbation frameworks of Gou\"ezel--Liverani~\cite{MR2201945} and Sedro--Rugh~\cite{SedroRugh2021}. Our emphasis is the separation 
of the highest order infinitesimal variation in the explicit response formulae, 
which permits a recursive construction of perturbations with prescribed responses,
as we shall explain below.

\medskip

This paper concentrates on direction (ii) in the category of 
smooth expanding endomorphisms on compact manifolds of arbitrary dimension, 
for which the physical measures are absolutely continuous invariant probability measures (ACIPs) with respect to volume.
Our first main result - Theorem~\ref{thm: response eend} -
 investigates the higher order differentiability of ACIP densities,
and derives explicit response formulae of any order. 
Roughly speaking,  
using the Riemannian exponential map, 
if a parameter family of smooth expanding endomorphisms $\{f_t\}$  
has a local expansion of the form
\beqn
f_t(x)=\exp_{f(x)}\left(t v_1(x)+t^2v_2(x)+\cdots+t^rv_r(x)
+o\left(|t|^r\right)\right) \ \ \text{as} \ t\to 0,
\eeqn
where  
$v_k$ are sections of the pull-back bundle $f^*TX$ for each $1\le k\le r$, 
we prove that the corresponding family of ACIP densities $\{\rho_t\}$ has a local expansion 
\beqn
\label{Taylor rhot0}
\rho_t=\rho+t\rho^{(1)}+t^2\rho^{(2)} + \dots + t^r\rho^{(r)} + o\left(|t|^r\right)
\ \ \text{as} \ t\to 0,
\eeqn
where $\rho^{(k)}$ is called the $k$-th order response at $f_0=f$.
Our explicit formula of $\rho^{(k)}$ matches with Ruelle's formal expression in~\cite{MR1492947},
but provides an alternative representation to emphasize its major term,
which only depends on $v_k$ but not on $v_1, \dots, v_{k-1}$. 

\medskip

The proof of Theorem~\ref{thm: response eend} relies on 
the cerebrated perturbation theory of Keller--Liverani and
Gou\"ezel--Liverani~\cite{MR1679080, MR2201945},
which ensures two crucial spectral properties~\ref{property P2} 
and \ref{property P3} for the family of transfer operators $\{\CL_t\}$. 
Together with Property~\ref{property P1} - weak differentiability of $\{\CL_t\}$ with respect to $t$, 
we are able to rigorously derive the expansion of the resolvent and leading eigen-projection of $\{\CL_t\}$.
See Section~\ref{sec: abs theory} for more details. 
 
\medskip 
 
In practice, it would be interesting to work on the inverse problem to the response formulae, that is, 
we seek a suitable perturbation to realize certain prescribed responses. 
More precisely,
given an expanding endomorphism $f$ and functions $\wrho_1, \dots, \wrho_r$, 
we intend to find a perturbation family $\{f_t\}$ of $f_0=f$ 
such that the $k$-th order response of 
the corresponding ACIP densities $\{\rho_t\}$
is exactly the given $\wrho_k$ for any $k=1, \dots, r$. 
Due to the explicit response formulae given by Theorem~\ref{thm: response eend},
it suffices to seek suitable infinitesimal variations $v_1, \dots, v_r$.  
Following the terminology in~\cite{MR3787356}, we call such inverse problem the ``request problem" of order $r$. 
If the request problem has solutions, we could further discuss the optimal solutions for some cost functionals. 

\medskip

In the class of circle expanding maps,
the linear request problem was well studied by Galatolo and Pollicott~\cite{MR3670005},
in which they proved the infinite-dimensionality of 
the solution space, and they also showed
the existence and uniqueness of minimal solution 
for the cost functional being $L^2$ or a suitable Hilbert Sobolev norm. 
Note that the existence of solutions to linear request  may not need hyperbolicity of the underlying dynamics, 
as shown by Kloeckner~\cite{MR3787356} using Moser's trick. 

\medskip

We also mention that there is another type  of optimization problem related to linear response,
that is, to find  infinitesimal variation $v_1$ to maximize the expectation $\ds \int c \rho^{(1)}$
for some observable $c: X\to \RR$. 
Recently, 
Froyland  and Galatolo~\cite{MR4863642} established the unique maximizing solution 
for the circle expanding maps,
and similar result was obtained by Antown, Froyland and Galatolo~\cite{MR4480803}
for certain integral operators and stochastic systems.

\medskip

Our second main result - Theorem~\ref{thm: request} gives an affirmative solution to 
the request problem of any order, in the category of smooth expanding endomorphisms on  compact manifolds. 
That is, 
given an unperturbed endomorphism $f$ and intended responses $\wrho_1, \dots, \wrho_r$, 
under some mild regularity assumption, 
solutions to the request problem (in terms of infinitesimal variations
$(v_1, \dots, v_r)$) form an infinite dimensional set. 
In fact, explicit expressions of all solutions satisfying
the additional regularity condition~\eqref{solution consistent} will be given by \eqref{explicit request sol}, by solving the equation system~\eqref{solve v} due to response formulae. 
As an immediate consequence, Corollary~\ref{cor: request opt} shows 
the existence of relaxed minimizers for some cost functional,
under suitable coercivity and lower semicontinuity assumptions.

Theorem~\ref{thm: request} and Corollary~\ref{cor: request opt} can be viewed as an extension of the results by
Galatolo and Pollicott~\cite{MR3670005} to higher order and to arbitrary  manifolds. 
To demonstrate the complication of higher order response and request, 
we shall explicitly compute the quadratic response formula 
for smooth circle expanding endomorphisms in Section~\ref{subsec: example general}, 
and we show in Section~\ref{subsec: opt quad doubling} that 
relaxed minimizers need not be unique.
We stress that new phenomena arise in the higher order request problem, 
since it involves nonlinear terms of infinitesimal variations (see \eqref{def w22} for instance). 
Also, Example~\ref{example: quadratic non-uniq} gives a perturbation invisible to linear response 
but with a nonzero second order response, even without a second order perturbation term.

\medskip

The paper is organized as follows. 
In Section~\ref{sec: 2}, we first introduce some basic notations and necessary
background for smooth expanding endomorphisms, 
then we state our two main results, that is, Theorem~\ref{thm: response eend} 
on the response formulae 
and Theorem~\ref{thm: request} on the request problem.
Section~\ref{sec: abs theory} 
is devoted to the perturbation theory of transfer operators,
including  expansions of the resolvent operator and the leading eigen-projection.
Next we prove Theorem~\ref{thm: response eend}  in Section~\ref{sec: proof1},
and prove Theorem~\ref{thm: request}  in Section~\ref{sec: proof2}.
Finally, we apply our main results to the quadratic response and request for circle expanding maps.

\section{Settings and Main Results}
\label{sec: 2}

\subsection{Basics for Expanding Endomorphisms}

Let $X$ be a compact connected orientable smooth Riemannian manifold  without boundary, 
endowed with a normalized Riemannian volume $m$. 
Given $\alpha> 1$,
a map $f: X\to X$ is said to be a \emph{$C^\alpha$ smooth  expanding endomorphism} of $X$,
if $f$ is $C^\alpha$ smooth
and
there exist $\Lambda>1$ and $C>0$ such that 
\beqn
\left| D_xf^n(v)\right|\ge C\Lambda^{n} |v|, \ \ \ \text{for any} \ n\in \NN \ \text{and any}\  (x, v)\in TX.
\eeqn
We  denote by   $\EEnd^\alpha(X)$
the space of all $C^\alpha$ smooth expanding endomorphisms of $X$.

\medskip

It is well known (see e.g.~\cite{MR2129258},\cite{MR1793194}) that any $f\in \EEnd^\alpha(X)$ admits a unique 
\emph{absolutely continuous invariant probability measure (ACIP)},
which has a strictly positive density function  $\rho\in C^{\alpha-1}(X)$.
In fact, the ACIP density 
$\rho$ is the fixed point of the \emph{transfer operator} $\CL$ of $f$ 
(with respect to $m$),
which is defined by the following duality:  
\beq
\label{def CL0}
\int    \CL\varphi(x)  \psi(x) dm(x)
=
\int    \varphi(x)  \psi(f(x)) dm(x),
\ \ \text{for any} \ \varphi\in L^1(m), \ \psi\in L^\infty(m).
\eeq
For any $\beta\in (0, \al-1]$, 
the action of $\CL$ restricted on $C^{\beta}(X)$ has \emph{spectral gap},
i.e., 
$\mathrm{Spectrum}(\CL)$ consists of a simple leading eigenvalue $1$
and a compact subset  strictly inside the unit disc.
Equivalently, by the classical perturbation theory of linear operators (see e.g.~\cite{MR1335452}),
we can write  
$\CL=\CP + \CN$
with the following properties:
\begin{enumerate}[(1)]
\item $\CP$ is a one-dimensional projection, i.e., $\CP^2=\CP$ and $\dim(\mathrm{Im}(\CP))=1$;
\item $\CN$ is a bounded linear operator with spectral radius strictly less than $1$;
\item $\CP\CN=\CN\CP=0$. 
\end{enumerate}
Furthermore, $\CP$ is the leading eigen-projection of $\CL$ given by the following formula:
\beq
\label{def CP0}
\CP\varphi =\rho \int \varphi \, dm,  
\ \ \text{for any} \ \
 \varphi\in C^{\beta}(X).
\eeq

\medskip

\subsection{Response Formulae of Any Order}

Let $r\ge 1$ and $\alpha>r+2$.
We assume that a family $\{f_t\}\subseteq \EEnd^{\alpha}(X)$ is $C^r$ smooth, 
where $t\in (-t_0, t_0)$ for some $t_0>0$, 
that is, $t\mapsto f_t$ is a $C^r$ smooth map from  $(-t_0,  t_0)$
to $\EEnd^{\alpha}(X)$. 
For convenience, we write $f=f_0$.
Since $X$ is a compact Riemannian manifold
without boundary, its exponential map $\exp:TX\to X$
is globally defined. 
Shrinking $t_0$ if necessary, we may 
assume that 
the geodesic distance between $f_t(x)$ and $f(x)$ is less than the injective radius of $X$
for all $(t, x)\in (-t_0, t_0)\times X$.
We therefore  define
\beqn
v(t, x):=\exp_{f(x)}^{-1}(f_t(x))\in T_{f(x)}X, \ 
\ \  \text{for any}  \
(t, x)\in (-t_0, t_0)\times  X,
\eeqn
where $\exp_{f(x)}^{-1}$ denotes the inverse of the
exponential map on the normal neighborhood of $f(x)$.
The mapping $t\mapsto v(t,\cdot)$ is a $C^r$ curve in the
Banach space $\Gamma^\alpha(f^*TX)$,
which consists of all $C^\alpha$ smooth sections
of the pullback bundle
\beqn
f^*TX
:=\{(x, w): \ x\in X,\ w\in T_{f(x)}X\}.
\eeqn
For any $1\le k\le r$,  we define
the \emph{(normalized) $k$-th  infinitesimal variation} $v_k$ of $\{f_t\}$ at $t=0$
by 
\beq
\label{def vk}
v_k(x):=
\left.\frac{1}{k!}\p_t^k 
 v(t, x) \right|_{t=0}
=
\left.\frac{1}{k!}\p_t^k
\left[\exp_{f(x)}^{-1}(f_t(x))\right]\right|_{t=0}\in T_{f(x)}X,
 \ \ \text{for any}\  x\in X.
\eeq
Here ``normalized'' refers to the division by $k!$.
We notice that $\ds v_k\in\Gamma^\alpha(f^*TX)$, i.e.,
$v_k$ is a $C^\alpha$ section of $\ds f^*TX$,
and thus we  rewrite
\beqn
\exp_{f(x)}^{-1}(f_t(x))=v(t, x)=
t v_1(x)+t^2v_2(x)+\cdots+t^rv_r(x)
+o(|t|^r) \ \ \text{as} \ t\to 0.
\eeqn
We stress that, when the Riemannian metric on $X$ is fixed,
the  infinitesimal variations
$v_k(x)$ are intrinsically
defined via the Riemannian exponential map at $f(x)$.
\medskip

Our first main result establishes the 
response formulae of any order
for the ACIP densities of smooth expanding endomorphisms.

\begin{theorem}
\label{thm: response eend}
Let $r\ge 1$ and $\alpha>r+2$.
Let $\{f_t\}$ be a $C^r$ smooth family in $\EEnd^{\alpha}(X)$, 
and let $\{\rho_t\}$ be the corresponding family of ACIP densities. 
For convenience, write $f=f_0$ and $\rho=\rho_0$.
Then  $t\mapsto\rho_t$ is a $C^r$ curve in $C^{\alpha-r-2}(X)$ near $t=0$, such that the normalized derivative
\beq
\label{explicit rhok}
\rho^{(k)}
:= \left.\frac{1}{k!}\p_t^k \rho_t \right|_{t=0}
=-\CR \dive\left( \CL\left(\rho v_k\right)\right) + \CR \rho^{(k)}_\# .
\eeq
belongs to $C^{\alpha-1-k}(X)$
for every $1\le k\le r$, where
\begin{enumerate}[(1)]
\item  
$\ds \CR:=\sum_{n=0}^\infty \CL^n(I-\CP)$,
in which $\CL$ is the transfer operator of $f$ given by \eqref{def CL0},
$I$ is the identity map, and 
$\CP$ is the leading eigen-projection of $\CL$ given by \eqref{def CP0}
(see~\eqref{def CR} for more details); 
\item 
$\dive(\cdot)$ is the divergence operator (see the precise definition in~\eqref{def divergence});
\item
$\rho^{(k)}_\#$  is determined by 
$f$, $v_1, \dots, v_{k-1}$, but it 
is independent of $v_k, v_{k+1}, \dots, v_r$. Moreover, 
$\rho^{(k)}_\#$ has vanishing integral, i.e., $\ds \int \rho^{(k)}_\# dm=0$.
\end{enumerate}
In particular, 
$\ds
\rho^{(1)}_\#\equiv 0
$
and thus
$\ds
\rho^{(1)}=-\CR\, \dive\left( \CL\left(\rho v_1\right)\right).
$
\end{theorem}

\begin{remark}
The proof of Theorem~\ref{thm: response eend} 
gives the following estimates: for each $1\le k\le r$,  
\[
\left\|\rho_t-\rho-\sum_{j=1}^k t^j\rho^{(j)}\right\|_{C^{\alpha-k-2}(X)}
=o(|t|^k).
\]
Nevertheless, 
the coefficients   retain the stronger regularity, i.e.,
$\rho^{(k)}\in C^{\alpha-k-1}(X)$.
\end{remark}

Formal expressions of higher order responses 
were derived by Ruelle in~\cite{MR1492947},
which is similar to that of \eqref{rhok2} below.
Theorem~\ref{thm: response eend} provides an alternative form  
with the emphasis on the major term of $\rho^{(k)}$, i.e., 
the term $ -\CR\, \dive\left( \CL\left(\rho v_k\right)\right)$.
This particular form plays a crucial role in solving the request problem - Problem~\ref{qu: request} in 
the next subsection. 
We stress that the remainder term $\rho^{(k)}_\#$ could be quite complicated 
since it may involve the curvature of the 
Levi-Civita connection of the compact manifold $X$, but it is very fortunate that we do not need its explicit formula   for the request problem.

\subsection{Optimal Request Problem}
 
The request problem is regarded as the inverse problem to the response formula. 
In the category of smooth expanding endomorphisms,
it can be formulated as follows:

\begin{problem}[Request Problem]
\label{qu: request}
Given   $f\in \EEnd^{\alpha}(X)$ and
 $\wrho^{(k)}\in C^{\alpha-1-k}(X)$ for $1\le k\le r$, is there a $C^r$ perturbation family 
$\ds \{f_t\}\subseteq\EEnd^{\alpha}(X)$, which admits  ACIP densities $\{\rho_t\}$, 
such that 
\beqn
f_0=f, \ \ \text{and} \ \ \ \ds \rho^{(k)}:=\frac{1}{k!}\left.\p_t^k \rho_t\right|_{t=0}=\wrho^{(k)}
\ \ \text{for any}\ 1\le k\le r ?
\eeqn
\end{problem}

\medskip

To this end, we introduce the notion of solution set to the request problem.

\begin{definition}
We say $\ds (v_1, \dots, v_r)\in\left[\Gamma^{\alpha}(f^*TX)\right]^r$ is a solution to 
Problem~\ref{qu: request} if the following holds:
for any $C^r$ smooth family $\ds \{f_t\}\subseteq \EEnd^{\alpha}(X)$ with $f_0=f$
such that its $k$-th infinitesimal variation  equals $v_k$ for $1\le k\le r$
(see \eqref{def vk}),
the corresponding ACIP densities $\{\rho_t\}$  have an expansion of the form
\beq
\label{rho Taylor}
\rho_t=\rho+t\wrho^{(1)}+t^2\wrho^{(2)} + \dots + t^r\wrho^{(r)} + o\left(t^r\right),
\ \text{as}\ t\to 0,
\eeq
where $\rho$ is the ACIP density of $f$.
The remainder in~\eqref{rho Taylor} is taken in $C^{\alpha-r-2}(X)$.
We denote 
the solution set of  Problem~\ref{qu: request}
by $\ds \CV\left(f; \wrho^{(1)}, \dots, \wrho^{(r)}\right)$.
\end{definition}

\medskip 

Our second main result establishes the existence of solutions to the request problem.

\begin{theorem}
\label{thm: request}
Let $r\ge 1$ and $\alpha>r+2$.
If $\dim X\ge2$, we also assume that $\alpha$ is not an integer. If the prescribed data
$\left\{f; \wrho^{(1)}, \dots, \wrho^{(r)} \right\}$ satisfies 
\begin{enumerate}
\item[(\textbf{R1})]
$\ds 
\int \wrho^{(k)}\, dm=0
$
for any  $1\le k\le r$;
\item[(\textbf{R2})]
$f\in \EEnd^{\alpha+r}(X)$, and
$\wrho^{(k)}\in C^{\alpha+r-1-k}(X)$ for any $1\le k\le r$,
\end{enumerate}
then the solution set $\ds \CV\left(f; \wrho^{(1)}, \dots, \wrho^{(r)}\right)$ is non-empty, in fact,  it is infinite dimensional.
\end{theorem}

 \medskip
 
Let us first remark on the assumptions of Theorem~\ref{thm: request}.
The noninteger restriction of $\alpha$ is used only for the elliptic regularity estimate when $\dim X\ge 2$ (see Remark~\ref{rem: noninteger} for details). As for the conditions on the prescribed data,
on the one hand,
the vanishing integral condition in Assumption \textbf{(R1)} is natural, since 
\beqn
\int \wrho^{(k)} dm = \frac{1}{k!}\left.\p_t^k \left(\int \rho_t dm \right) \right|_{t=0}  = 
 \frac{1}{k!} \p_t^k 1=0
\eeqn
if \eqref{rho Taylor} holds.
On the other hand,
in the statements of Theorem~\ref{thm: response eend},
the unperturbed endomorphism $f$ is assumed to be only $C^\alpha$ smooth
and the $k$-th order response $\rho^{(k)}$ is supposed to be only $C^{\alpha-1-k}$ smooth,
yet we raise the smoothness by $r$ for both $f$ and $\wrho^{(k)}$ 
in  Assumption \textbf{(R2)}.
The reason to impose higher regularity on   $\left\{f; \wrho^{(1)}, \dots, \wrho^{(r)} \right\}$
may not be obvious at first sight, 
but we shall see in the proof that
it ensures that the recursively constructed variations belong to $\left[\Gamma^{\alpha}(f^*TX)\right]^r$
(see \eqref{solution consistent}).

\medskip

We then make a comment on the conclusion of Theorem~\ref{thm: request}.
Albeit $\ds \CV\left(f; \wrho^{(1)}, \dots, \wrho^{(r)}\right)$ need not be a linear subspace of $\left[\Gamma^{\alpha}(f^*TX)\right]^r$,
the infinite dimensionality means that it contains an infinite dimensional
affine subspace of $[\Gamma^\alpha(f^*TX)]^r$.
We also remark that,
as a matter of fact, all solutions in
$\ds \CV\left(f; \wrho^{(1)}, \dots, \wrho^{(r)}\right)$
satisfying~\eqref{solution consistent}
can be explicitly written down,
see \eqref{explicit request sol} for its precise formulation.

\medskip

The following corollary shows that solutions to the relaxed request
problem can be chosen optimally with respect to a suitable cost functional.

\begin{cor}
\label{cor: request opt}
Under the same assumptions of Theorem~\ref{thm: request}, 
let $\CV:=\CV\left(f;\wrho^{(1)},\dots,\wrho^{(r)}\right)$ and $\overline{\CV}^{\,\beta}$ be the closure of $\CV$ in $[\Gamma^\beta(f^*TX)]^r$ where $\beta\in(r+1,\alpha)$. 
Suppose that a function $\CJ:[\Gamma^\alpha(f^*TX)]^r\to\RR$ satisfies the following conditions:
	\begin{enumerate}[(1)]
		\item
		$\CJ$ is bounded below on $\CV$ and coercive with respect to
the $C^\alpha$-norm,
		that is,
		\beqn
		\inf_{\bv\in\CV}\CJ(\bv)>-\infty\ \text{ and }
		\lim_{\substack{\bv\in\CV\\ \|\bv\|_\alpha\to\infty}}
		\CJ(\bv)=+\infty,
		\eeqn
		where $\|\bv\|_\alpha:=\|v_1\|_\alpha+\dots+\|v_r\|_\alpha$ for $\bv=(v_1,\dots,v_r)$.
		\item
		$\CJ$ is lower semicontinuous on $\CV$ in the topology induced by $[\Gamma^\beta(f^*TX)]^r$, that is, for every $\bv_n,\bv\in\CV$ with $\|\bv_n-\bv\|_\beta\to0$, we have
		$\ds
		\CJ(\bv)\le\liminf_{n\to\infty}\CJ(\bv_n).
		$
	\end{enumerate}
Let $\overline{\CJ}:\overline{\CV}^{\,\beta}\to\RR\cup\{+\infty\}$ be the extension of $\CJ$, defined by
	\beqn
	\overline{\CJ}(\bv)
	:=\inf_{\substack{\bv_n\in\CV\\ \|\bv_n-\bv\|_\beta\to0}}
	\liminf_{n\to\infty}\CJ(\bv_n).
	\eeqn
Then there exists $\overline{\bv}=(\ov_1,\dots,\ov_r)\in\overline{\CV}^{\,\beta}$ such that
	\beqn
	\overline{\CJ}(\overline{\bv})
	=\min\left\{\overline{\CJ}(\bv): \  \bv\in\overline{\CV}^{\,\beta}\right\}
	=\inf\left\{\CJ(\bv):\ \bv\in\CV\right\}.
	\eeqn
\end{cor}

Corollary~\ref{cor: request opt} implies
the existence of a minimizer for the extended cost functional 
$\overline{\CJ}$ in
$\overline{\CV}^{\,\beta}$,
while the uniqueness of such a minimizer
is not guaranteed in general,
except in some particular situations.
For the doubling map, Galatolo and Pollicott~\cite{MR3670005}
obtained a unique solution of minimum $L^2$ norm to the linear
request problem. For general smooth expanding circle maps, they
established uniqueness for a suitable Hilbert Sobolev norm.
In Subsection~\ref{subsec: opt quad doubling}, 
we shall investigate the quadratic request problem in the same category 
and establish explicit formulae for the solution set 
$\CV\left(f; \wrho^{(1)}, \wrho^{(2)}\right)$;
we also provide an example for which 
the uniqueness of minimal solution could fail (see Example~\ref{example: quadratic non-uniq}).

\medskip 
 
Let us remark on the proof of Theorem~\ref{thm: request}
and Corollary~\ref{cor: request opt}.
To construct solutions, we 
convert the request problem to the equation system~\eqref{solve v1},
of which the triangular structure  allows us to solve \eqref{def k-th eq} recursively. 
For the optimization problem,
one could make a finite-dimensional approximation of the exact solutions
given in~\eqref{explicit request sol}, over which one may minimize.
Convergence of these finite-dimensional approximations to a
minimizer requires a separate argument.

We also remark that
the topology of $X$ alone does not quite determine uniqueness of an optimal solution.
For the doubling map, the linear request problem has a unique minimizer for $L^2$-norm functional, whereas Example~\ref{example: quadratic non-uniq} gives distinct relaxed minimizers for a quadratic request problem. Uniqueness therefore depends also on the cost functional and the request constraints.


\section{Perturbation Theory of Transfer Operators}
\label{sec: abs theory} 

\subsection{Properties  of Transfer Operators}

Let $r\ge 1$ and $\alpha>r+2$, and
let $\{f_t\}$ be a $C^r$ smooth family in $\EEnd^{\alpha}(X)$, 
where $t\in (-t_0, t_0)$ for some $t_0>0$.
We denote by
$\{\CL_t\}$ the corresponding family of transfer operators, 
which is characterized by the following duality relation
with respect to the Riemannian volume $m$:
\beq
\label{eq: duality}
\int    \CL_t \varphi(x)  \psi(x) dm(x)
=
\int    \varphi(x)  \psi (f_t(x)) dm(x),
\eeq
for any $t\in (-t_0, t_0)$, $\varphi\in L^1(m)$ and $\psi\in L^\infty(m)$.
Since each $f_t$ is locally diffeomorphic,  
decreasing the value of $t_0$ if necessary, 
we may assume that
all $f_t$ have the same number of local inverse branches, 
say $d$ branches, and 
write $f_t^{-1}(x):=\{y_t^{\ell}(x)\}_{1\le \ell\le d}$ for $t\in (-t_0,t_0)$ and  $x\in X$.
Moreover, 
$(t, x)\mapsto y_t^{\ell}(x)$ is $C^{r}$ smooth in $t$ and $C^{\alpha}$ smooth in $x$.  
By change of variables in \eqref{eq: duality},
the transfer operator $\CL_t$ can be rewritten as
\beq
\label{def transfer2}
\CL_t \varphi(x)=\sum_{\ell=1}^d  \varphi(y_t^{\ell}(x))\left| \det\left(D y_t^{\ell}(x)\right)\right|,
\ \ \text{for any} \ \varphi\in L^1(m),
\eeq
where $\ds D y_t^{\ell}(\cdot)$ denotes the tangent map of $y_t^{\ell}(\cdot)$.

 \medskip

We shall restrict the action of $\CL_t$ on finer subspaces of $L^1(m)$.
More precisely, we consider a nested sequence of function spaces $\{\CB_i\}_{-1\le i\le r}$, where
$\CB_i:=C^{\alpha-r-1+i}(X)$.
In particular, $\CB_{-1}=C^{\alpha-r-2}(X)$ has positive regularity.
It immediately follows from \eqref{def transfer2} that 
$\CL_t(\CB_i)\subseteq \CB_i$ and the restriction  $\CL_t|_{\CB_i}$ is a bounded linear operator
for any $t\in (-t_0, t_0)$.
Moreover, the family $\ds \{\CL_t\}$ of transfer operators satisfies
the following property:

\medskip

\begin{propertyx}
\label{property P1}
\textbf{Weak differentiability at $t=0$}:
there exist   $\CL^{(1)}, \CL^{(2)}, \dots, \CL^{(r)}$ such that the following statements hold
for any 
$1\le k\le  i \le r$:
\begin{enumerate}[(1)]
\item $\CL^{(k)}$ is   a bounded linear operator from $\CB_i$ to $\CB_{i-k}$;
\item $ \left\|\CL_t -\CL - \sum_{j=1}^k t^j \CL^{(j)}\right\|_{\CB_i\to \CB_{i-k-1}}=o\left( t^{k}\right)$,
as $t\to 0$.
\end{enumerate}
Moreover, for every $0\le i\le r$, we have
\(\ds
\|\CL_t-\CL\|_{\CB_i\to\CB_{i-1}}=O(|t|)
\)
as $t\to 0$.
\end{propertyx}

\begin{remark}
Property~\ref{property P1} follows by differentiating
\eqref{def transfer2} with respect to $t$ and applying Taylor's
formula, together with compact embeddings and H\"older interpolation
after an additional loss of spatial regularity.
In fact, any additional loss $\epsilon\in(0,1]$ suffices
for~\ref{property P1}(2): the remainder is $o(|t|^k)$
in the operator norm from $\CB_i$ to
$C^{\alpha-r-1+i-k-\epsilon}(X)$.
See also Sedro--Rugh~\cite{SedroRugh2021}
on a little H\"older scale.
Here we take $\epsilon=1$, which is sufficient for our applications.
\end{remark}

Explicit formulae for  $\CL^{(1)}, \CL^{(2)}, \dots, \CL^{(r)}$ will be given by
Lemma~\ref{lem: explicit CLk}.
Note that 
$\CL_t$ compatibly acts on  
$\{\CB_i\}_{-1\le i\le r}$, i.e.,
$\left(\CL_t|_{\CB_j}\right)|_{\CB_i}=\CL_t|_{\CB_i}$
for any $-1\le j\le i\le r$.
Therefore, we simply write $\CL_t$ instead of $\CL_t|_{\CB_i}$
when $\CB_i$ is clear.

\medskip

In the rest of this subsection, 
we shall  introduce two spectral properties -~\ref{property P2} and~\ref{property P3} 
below,
which are compatible on $\{\CB_i\}_{-1\le i\le r}$,
and thus we omit the underlying space.
After shrinking \(t_0\) if necessary,
these properties follow from uniform Lasota--Yorke estimates and weak continuity, together with the
perturbation theory 
due to Keller-Liverani~\cite{MR1679080}
and Gou\"ezel-Liverani~\cite{MR2201945}, 
see also Appendix A.3 in Baladi's book~\cite{MR3837132}.

\medskip

\begin{propertyx}
\label{property P2}
\textbf{Uniform spectral gap}: 
there exists $\eta\in (0, 1)$ such that 
for any $t\in (-t_0, t_0)$, 
$\mathrm{Spectrum}(\CL_t)$ consists of a simple leading eigenvalue $1$
and a compact subset  inside the disc of radius $\eta$.
Equivalently,  
we can  write  
$\CL_t=\CP_t + \CN_t$
with the following properties:
\begin{enumerate}[(1)]
\item $\CP_t$ is a one-dimensional projection, i.e., $\CP_t^2=\CP_t$ and $\dim(\mathrm{Im}(\CP_t))=1$;
\item $\CN_t$ is a bounded linear operator with spectral radius less than $\eta$;
\item $\CP_t\CN_t=\CN_t\CP_t=0$. 
\end{enumerate}
\end{propertyx}

\medskip

Let $I$ denote the identity map, 
and recall that the resolvent of $\CL_t$ is defined by
\beqn
R_t(z):=\left(z I -\CL_t\right)^{-1},
\ \ \text{whenever} \ \
z\in \CC\backslash \mathrm{Spectrum}(\CL_t).
\eeqn 
We stress that for any fixed $z$, 
the mapping $t\mapsto R_t(z)$ might not be continuous when 
$R_t(z)$ is regarded as an operator from $\CB_i$ to itself.
Nevertheless,  we have

 \medskip
 
\begin{propertyx}
\label{property P3}
\textbf{Uniform boundedness of the resolvent along a circle enclosing  $1$}: 
there exists $\ds \delta\in (0, (1-\eta)/2)$ 
such that 
\beq
\label{eq: Rt bounded}
\sup_{-1\le i\le r}\sup_{t\in (-t_0, t_0)} \sup_{z\in \gamma} \|R_t(z)\|_{\CB_i\to \CB_i} <\infty,
\eeq
where  
$\ds 
\gamma:=\left\{ z\in \CC: |z-1|=\delta\right\}
$.
 
\end{propertyx} 
 
\medskip

In fact, a stronger uniform boundedness 
with $\gamma$ replaced by $V_{\delta}$ in~\eqref{eq: Rt bounded} 
was obtained
by Gou\"ezel-Liverani~\cite{MR2201945}, 
where  
$\ds 
V_{\delta}:=\left\{ z\in \CC: |z|\ge \eta+\delta, \ 
\text{and} \ |z-1|\ge \delta\right\}
$
for some $\delta\in (0, (1-\eta)/2)$.
Anyhow, for our purpose, we only require a weaker condition \ref{property P3}.

\medskip

In the next subsections, we shall 
derive the expansion of the resolvent $R_t(z)$ (see Lemma~\ref{lem: der Rt}),
and the expansion of the leading eigen-projection $\CP_t$ (see Lemma~\ref{lem: der CP}),
with respect to $t$ near $t=0$.
We remark that the proof of these two lemmas
actually works for any operators satisfying Properties~\ref{property P1}-\ref{property P3},
which need not be the transfer operators.

 
\subsection{Expansion of the Resolvent}

Let $\delta$ be given by Property  \ref{property P3},
and set
$U:=\{z\in \CC:  |z-1|< 2\delta\}$.
It follows from Property  \ref{property P2} that
$\ds \mathrm{Spectrum}(\CL_t)\cap U=\{1\}$.
The following lemma shows that 
the resolvent $R_t(z)$ of $\CL_t$ is 
in fact a meromorphic function on $U$ with a simple pole at $z=1$.

\begin{lemma}
For any $t\in (-t_0, t_0)$ and 
any $z\in U\backslash\{1\}$, 
we have
\beq
\label{def exp Rt}
R_t(z)=\sum_{n=-1}^\infty (z-1)^n\CR_t^{(n)},
\ \ \text{where} \ \ 
\CR_t^{(n)}=
\begin{cases}
\CP_t, & n=-1, \\
(-1)^n(I-\CN_t)^{-n-1}(I-\CP_t), & n\ge 0.
\end{cases}
\eeq
\end{lemma}

\begin{proof}
Property \ref{property P2} implies that
$\ds
\CL_t\CP_t=\CP_t\CL_t= \CP_t
$
and
$\ds
\CL_t(I-\CP_t)=(I-\CP_t)\CL_t=\CN_t.
$
Therefore, $R_t(z)=R_t(z)\CP_t +  R_t(z)(I-\CP_t)$ for any $z\in U\backslash\{1\}$, where
\beqn
\begin{split}
& R_t(z)\CP_t=(z I-\CL_t)^{-1} \CP_t =(z-1)^{-1} \CP_t, 
\\
& R_t(z)(I-\CP_t)=(z I-\CL_t)^{-1} (I-\CP_t)=(z I-\CN_t)^{-1} (I-\CP_t).
\end{split}
\eeqn
By the spectral mapping theorem (see e.g. Chapter VII in \cite{MR1070713}), we have
$$\ds \mathrm{Spectrum}((I-\CN_t)^{-1})
=\left(1-\left(\mathrm{Spectrum}(\CN_t)\right)\right)^{-1},
$$
and thus the spectral radius of 
$(z-1)(I-\CN_t)^{-1}$ is less than $2\delta/(1-\eta)<1$.
Then we have
\beqn
\begin{split}
(zI-\CN_t)^{-1} (I-\CP_t)
& =(I-\CN_t)^{-1} \left[ I +(z-1)(I-\CN_t)^{-1}\right]^{-1} (I-\CP_t) \\
&=\sum_{n=0}^\infty (-1)^n(z-1)^n (I-\CN_t)^{-n-1}(I-\CP_t).
\end{split}
\eeqn
Therefore, 
\beqn
R_t(z)=(z-1)^{-1} \CP_t + \sum_{n=0}^\infty (-1)^n(z-1)^n (I-\CN_t)^{-n-1}(I-\CP_t)=
\sum_{n=-1}^\infty (z-1)^n \CR_t^{(n)}.
\eeqn
The proof of this lemma is complete.
\end{proof}

\medskip

Let $\gamma$ be the circle given by Property \ref{property P3},
and note that $\gamma\subset U\backslash\{1\}$.
To obtain the expansion of
the mapping $t\mapsto R_t(z)$ at $t=0$ for any $z\in \gamma$, 
we introduce the following notations.
For any positive integer $k\in \NN$, we 
denote by $\Par(k)$ the set of all \emph{ordered partitions} of $k$, i.e.,
\beqn
\Par(k):=\left\{ (i_1  \dots i_\ell):
\ \ell, i_1, \dots, i_\ell\in \NN, 
 \ \text{and} \ 
\ i_1+\dots +i_\ell=k \right\}. 
\eeqn
Let $k$ and $\ell$ be two positive integers such that $1\le \ell\le k\le r$,
and recall that $\CL^{(i)}$ are given by Property \ref{property P1} for $1\le i\le r$.
For any ordered partition  $\bi=(i_1 i_2 \dots i_\ell)\in \Par(k)$, 
we associate an $(\ell+1)$-multilinear form defined by
\beq
\label{def CE}
\CE^{\bi}\left( \Phi_1, \Phi_2, \dots, \Phi_{\ell}, \Phi_{\ell+1}\right):=
\Phi_{1} \CL^{(i_1)} \Phi_{2}\CL^{(i_2)}  \dots \Phi_{\ell} \CL^{(i_\ell)} \Phi_{\ell+1}.
\eeq
Write $R(z)=R_0(z)$. For any $z\in \gamma$, we further define
\beq
\label{def Rkz}
R^{(k)}(z):=\sum_{\bi \in \Par(k)} 
\CE^{\bi}\left(  R(z),\dots,R(z)\right).
\eeq
 
\medskip

\begin{lemma}
\label{lem: der Rt}
For   any $z\in \gamma$ and any 
$1\le k\le i\le r$,
we have 
\beqn
\left\| R _t(z) - R(z) - \sum_{j=1}^k t^j R^{(j)}(z)\right\|_{\CB_i\to \CB_{i-k-1}}
=o\left( t^k \right), \
\text{as} \ t\to 0,
\eeqn
where $R^{(j)}(z)$ are given by \eqref{def Rkz}.
The remainder is uniform for $z\in\gamma$.
\end{lemma}

\begin{proof}
Set $R^{(0)}(z)=R(z)$,
$E_0(t)=\CL_t-\CL$ and
$\ds E_\ell(t)=\CL_t-\CL-\sum_{j=1}^\ell t^j\CL^{(j)}$ for $\ell\ge1$.
By~\eqref{def Rkz}, we have
\(\ds
R^{(n)}(z)=R(z)\sum_{j=1}^n\CL^{(j)}R^{(n-j)}(z),
\)
and 
thus $R^{(n)}(z):\CB_i\to\CB_{i-n}$ is  uniformly bounded on
$\gamma$ whenever $i\ge n$. For
$\ds S_k(t,z)=\sum_{n=0}^kt^nR^{(n)}(z)$, this recursion gives the exact identity
\beq
\label{eq: res Sk}
I-(zI-\CL_t)S_k(t,z)
=\sum_{n=0}^kt^nE_{k-n}(t)R^{(n)}(z).
\eeq
For any $n<k$, Property~\ref{property P1} gives
\[
\|E_{k-n}(t)\|_{\CB_{i-n}\to\CB_{i-k-1}}=o(|t|^{k-n})
\ \ \text{and} \ \
\|E_0(t)\|_{\CB_{i-k}\to\CB_{i-k-1}}=O(|t|).
\]
Consequently the right-hand side of \eqref{eq: res Sk}
is $o(|t|^k)$ from $\CB_i$ to
$\CB_{i-k-1}$, uniformly on $\gamma$. 
Multiplication from the left on both sides of \eqref{eq: res Sk} by
$R_t(z)$, which is bounded on $\CB_{i-k-1}$ due to Property~\ref{property P3},
proves the lemma.
\end{proof}

\subsection{Expansion of the Leading Eigen-projection}

By the theorem of ``separation of spectrum'' (see e.g. Chapter Three, \S 6.4 in~\cite{MR1335452}),
the leading eigen-projection $\CP_t$ 
given by Property \ref{property P2}
has the following representation: 
\beq
\label{CPt integral}
\CP_t =\dfrac{1}{2\pi i} \oint_\gamma R_t(z) dz,
\eeq
where $\gamma$ is the circle given by
Property \ref{property P3}.
For convenience, 
we shall write 
$\CP=\CP_0$, $\CN=\CN_0$ and $R(z)=R_0(z)$.
Taking $t=0$ in Formula~\eqref{def exp Rt},
we obtain the Laurent series of $R(z)$ at $z=1$:
\beq
\label{def exp R0}
R(z)=\sum_{n=-1}^\infty (z-1)^n \CR^{(n)},
\ \ \text{where} \ \ 
\CR^{(n)}:=
\begin{cases}
\CP, & n=-1, \\
(-1)^n(I-\CN)^{-n-1}(I-\CP), & n\ge 0.
\end{cases}
\eeq
Note that $\CP \CR^{(n)} = \CR^{(n)} \CP=0$ for any $n\ge 0$.

\medskip

To derive the expansion of
the mapping $t\mapsto \CP_t$ at $t=0$, 
we need the following notation.
Denote by $\ZZ_{\ge -1}$ the set of all integers no less than $-1$.
For any $\ell\ge 1$ and $s\in \ZZ$,
we define
\beqn
\Omega_\ell^s:=
\left\{ (j_1 j_2 \dots  j_{\ell}):
\  j_1, j_2, \dots, j_{\ell}\in \ZZ_{\ge -1}, 
 \ \text{and} \ 
\ j_1+j_2+\dots +j_{\ell}=s \right\}. 
\eeqn
That is, $\Omega_\ell^s$ is the set of all $\ell$-tuples of integers in $\ZZ_{\ge -1}$
with sum equal to $s$.
Given any
$\bi=(i_1 i_2 \dots i_\ell)\in \Par(k)$ with length $|\bi|=\ell$
and any
$\bj=(j_1 j_2 \dots  j_{\ell+1})\in \Omega_{\ell+1}^s$,
we define 
\beq
\label{def CRt ell}
\CE^{\bi, \bj}:=\CE^\bi(\CR^{(j_1)}, \CR^{(j_2)}, \dots, \CR^{(j_{\ell+1})} )
=\CR^{(j_1)} \CL^{(i_1)} \CR^{(j_2)} \CL^{(i_2)}  \dots \CR^{(j_\ell)}  \CL^{(i_\ell)} \CR^{(j_{\ell+1})},
\eeq
where $\CE^\bi(\cdots)$ is the $(\ell+1)$-multilinear form given by \eqref{def CE}.
We further define
\beq
\label{der Pt k0}
\CP^{(k)}= \sum_{\bi \in \Par(k)} 
\sum_{\bj\in \Omega_{|\bi|+1}^{-1}} \, \CE^{\bi, \bj},
\eeq
which is a bounded linear operator from $\CB_i$ to $\CB_{i-k}$
for any $1\le k\le i \le r$.

\medskip

The following lemma provides 
the expansion of 
the leading eigen-projection $\CP_t$ at $t=0$.

\begin{lemma}
\label{lem: der CP}
For  any $1\le k\le i\le r$,
we have 
\beqn
\left\| \CP_t - \CP - \sum_{j=1}^k t^j \CP^{(j)}\right\|_{\CB_i\to \CB_{i-k-1}}
=o\left( t^k \right), \
\text{as} \ t\to 0,
\eeqn
where $\CP^{(j)}$ are given by \eqref{der Pt k0}.
\end{lemma}

\begin{proof}
By~\eqref{def Rkz}, ~\eqref{CPt integral} and Lemma~\ref{lem: der Rt}, 
we have
\beqn
\begin{split}
\CP_t-\CP 
=\dfrac{1}{2\pi i} \oint_\ga \left[R_t (z) - R(z)\right] dz 
&=\sum_{j=1}^k t^j \ \dfrac{1}{2\pi i} \oint_\ga  R^{(j)} (z)   dz + o\left( t^{k}\right)\\
&=
\sum_{j=1}^k t^j \
\sum_{\bi \in \Par(j)} 
\dfrac{1}{2\pi i} \oint_\ga \CE^{\bi}\left( R(z),  \dots, R(z)\right) dz + o\left( t^{k}\right),
\end{split}
\eeqn
The remainder is taken from $\CB_i$ to $\CB_{i-k-1}$,
by integrating the estimate of Lemma~\ref{lem: der Rt}.
By~\eqref{def CE} and~\eqref{def exp R0},  
\beqn
\CE^{\bi}\left( R(z),  \dots, R(z)\right)
=\CE^{\bi}\left( \sum_{n=-1}^\infty (z-1)^n \CR^{(n)} ,  \dots, 
\sum_{n=-1}^\infty (z-1)^n \CR^{(n)} \right) 
\eeqn
is a meromorphic function on $U=\{z\in \CC:  |z-1|<2\delta\}$ with a possible pole at $z=1$ of order at most $(1+|\bi|)$,
and it is easy to check that its order one term is given by
$\ds
\dfrac{\sum_{\bj\in \Omega_{|\bi|+1}^{-1}} \, \CE^{\bi, \bj}}{z-1},
$
where $\ds \CE^{\bi, \bj}$
is defined by \eqref{def CRt ell}.
By the residue theorem, 
\beqn
\sum_{\bi \in \Par(j)}  \dfrac{1}{2\pi i} \oint_\ga \CE^{\bi}\left( R(z),  \dots, R(z)\right) dz
=
\sum_{\bi \in \Par(j)}  \sum_{\bj\in \Omega_{|\bi|+1}^{-1}} \, \CE^{\bi, \bj},
\eeqn
from which the lemma immediately follows.
\end{proof}


\section{Proof of Response Formulae of Any Order }
\label{sec: proof1}

\subsection{\texorpdfstring{Simplification of $\CP^{(k)}$}{Simplification of the projection coefficients}}

We recall the precise settings of Theorem~\ref{thm: response eend}, 
that is,
$\{f_t\}$ is a $C^r$ smooth family in $\EEnd^{\alpha}(X)$,
$\{\rho_t\}$ is the corresponding family of ACIP densities,
$\{\CL_t\}$ is the family of transfer operators 
corresponding to $\{f_t\}$,
and 
$\{\CP_t\}$ is the family of leading eigen-projections of $\{\CL_t\}$.
For convenience, we write 
$f=f_0$, $\rho=\rho_0$, $\CL=\CL_0$ and $\CP=\CP_0$. 

\medskip

The  resolvent $R(z)=(zI-\CL)^{-1}$  
is meromorphic on 
$U:=\{z\in \CC:  |z-1|<2\delta\}$ with a simple pole $z=1$, 
and its Laurent series expansion at $z=1$ is given by~\eqref{def exp R0}. 
Note that  the order zero coefficient of~\eqref{def exp R0} can be simplified as follows:
\beq
\label{def CR}
\CR:=\CR^{(0)}=(I-\CN)^{-1}(I-\CP)
=\sum_{n=0}^\infty \CN^n (I-\CP)
=\sum_{n=0}^\infty \CL^n (I-\CP).
\eeq

The expansion of leading eigen-projection
$\CP_t$ at $t=0$ has been established
in Lemma~\ref{lem: der CP},
with the coefficients $\CP^{(k)}$ given by Expression~\eqref{der Pt k0}.
We stress that~\eqref{der Pt k0}
works for any abstract operator family satisfying  Properties~\ref{property P1}-\ref{property P3},
while further simplification can be made in the case of transfer operators.
To be precise, we have the following lemma.

\begin{lemma}
\label{lem: CPkP}
Let $\ds \CR$ be given by \eqref{def CR}.
For  any $1\le k \le i \le r$, we have
\beq
\label{der Pt k}
\CP^{(k)}= \sum_{(i_1 i_2\dots i_\ell) \in \Par(k)} 
\CR \CL^{(i_1)}\CR \CL^{(i_2)} \dots \CR \CL^{(i_\ell)}  \CP,
\eeq
which is a bounded linear operator from $\CB_i$ to $\CB_{i-k}$.
\end{lemma}

\begin{proof}
Taking  $\psi\equiv 1$ in \eqref{eq: duality}, we obtain 
$\ds
\label{eq: dual0}
\int \CL_t \varphi \,  dm =\int \varphi \,  dm
$
for any
$\varphi\in L^1(m)$.
It then follows from Property \ref{property P1} that for
any $1\le j\le k\le i\le r$, we have:
\beq
\label{CLj vanishing}
\int \CL^{(j)}\varphi\, dm = 0, 
\ \ \text{for any} \ \
\varphi\in \CB_i:=C^{\al-r-1+i}(X).
\eeq
Combining \eqref{CLj vanishing} and \eqref{def CP0}, 
we get
$\ds \CP\CL^{(j)}\varphi =  \rho \int \CL^{(j)}\varphi\, dm=0$,
that is,
$\CP\CL^{(j)}\equiv 0$. 
Due to this fact
and the expressions of $\CR^{(j)}$ in \eqref{def exp R0},
for any $\bi=(i_1 i_2\dots i_\ell) \in \Par(k)$
and for any $\bj=(j_1 j_2 \dots j_\ell j_{\ell+1})\in \Omega_{|\bi|+1}^{-1}$, 
the term 
\beqn
\CE^{\bi, \bj}=
\CR^{(j_1)} \CL^{(i_1)} \CR^{(j_2)} \CL^{(i_2)}  \dots \CR^{(j_\ell)}  \CL^{(i_\ell)} \CR^{(j_{\ell+1})}
\eeqn
must be zero unless $\bj=(0, 0, \dots 0, -1)$. 
That is, the only possible non-zero term is 
\beqn
\CR^{(0)} \CL^{(i_1)} \CR^{(0)} \CL^{(i_2)}  \dots \CR^{(0)}  \CL^{(i_\ell)} \CR^{(-1)}
=
\CR \CL^{(i_1)}\CR \CL^{(i_2)} \dots \CR \CL^{(i_\ell)}  \CP.
\eeqn
Hence 
the expression of $\CP^{(k)}$ in  \eqref{der Pt k0} is simplified to \eqref{der Pt k}.
\end{proof}

\subsection{\texorpdfstring{Major Term of $\CL^{(k)}$}{Major Term of the transfer coefficients}}

Recall Property~\ref{property P1} that
$\CL_t$ is $C^r$ weakly differentiable at $t=0$.
The $k$-th operator derivative is taken from $\CB_i$ to $\CB_{i-k-1}$; 
its coefficient $\CL^{(k)}$ is also bounded into $\CB_{i-k}$.
We shall provide an explicit formula of 
$\CL^{(k)}=\frac{1}{k!}\left. \p_t^k\CL_t\right|_{t=0}$ for all $1\le k\le r$,
which is split into two parts: 
the first part only involves the $k$-th infinitesimal variation $v_k$ (see the definition in \eqref{def vk}),
and the second part only depends on  infinitesimal variations of order less than $k$.
Moreover, there is only one term in the first part,
which we call the major term of $\CL^{(k)}$.

To this end, we recall some definitions and notations in Riemannian geometry (see e.g.~\cite{MR2243772} for reference).
We denote the Riemannian metric on the smooth manifold $X$ by $\lan\cdot, \cdot\ran$, 
that is, $ \lan\cdot, \cdot\ran_x$ is an inner product on the tangent space $T_x X$ for each $x\in X$ 
such that $x\mapsto \lan\cdot, \cdot\ran_x$ is smooth.
Let $dm$ denote the smooth Riemannian volume form.
Given a  $C^1$ smooth vector field $v\in \Gamma^1(TX)$,  
the divergence of $v$ 
is defined to be the scalar $\dive(v)$ such that
\beq
\label{def divergence}
L_v dm = \dive(v) dm,
\eeq
where $L_v$ denotes the Lie derivative along the vector field $v$.
Since $X$ is a compact orientable smooth  manifold without boundary,
the following properties hold for the divergence:
\begin{enumerate}[(1)]
\item 
The divergence theorem implies that
\beq
\label{dive Lib0}
\int \dive(v) dm=0,
\ \ \text{for any} \ \ v\in \Gamma^1(TX).
\eeq
More generally, since 
$\ds
\ds \dive(\varphi v)=\varphi \dive(v) + \lan v, \nabla \varphi\ran
$
for any $\varphi\in C^1(X)$ and any $v\in \Gamma^1(TX)$, 
where $\nabla \varphi$ denotes the gradient vector field of $\varphi$,
again
by divergence theorem, we get
\beq
\label{dive Lib}
\int \lan v, \nabla \varphi\ran \, dm = -\int \varphi\, \dive(v)\, dm;
\eeq

\item If a real-valued function $\varphi\in C^\beta(X)$  satisfies $\ds \int \varphi dm=0$, 
where $\beta>0$ is not an integer (or $\beta\ge 0$ if $\dim X=1$),
then there exists a vector field $v\in \Gamma^{\beta+1}(TX)$ such that 
\beq
\label{div eq}
\dive(v)=\varphi.
\eeq
Indeed, let $d=\dim X$, by de Rham cohomology, 
the $d$-form $\varphi dm$ has vanishing integral if and only if it is exact, 
i.e., 
there is a $(d-1)$-form $\psi$ such that $d\psi =\varphi dm$.
Moreover, there is a vector field $v$ such that 
$\iota_v dm = \psi$, where $\iota_v$ denotes the contraction  by $v$.
Then \eqref{div eq} follows from Cartan's magic formula, 
since
\beqn
\dive(v) dm=L_v dm= (d\iota_v + \iota_v d) dm= d(\iota_v dm)=d\psi= \varphi dm.
\eeqn
\end{enumerate}

\begin{remark}
\label{rem: noninteger}
For the asserted regularity of $v\in \Gamma^{\beta+1}(TX)$ in Equation~\eqref{div eq}, 
a more direct method is to solve the Poisson equation $\Delta u=\varphi$.
By \cite[Theorem~4.7]{Aubin1998}, it admits a unique zero-mean solution
$u\in C^{\beta+2}(X)$, and hence $v=\nabla u$ belongs to
$\Gamma^{\beta+1}(TX)$ and satisfies $\dive (v)=\varphi$.
The Schauder estimate used here applies to noninteger
$\beta>0$ when $\dim X\ge2$.
In dimension one, we can drop such noninteger restriction, since
periodic integration in arclength coordinates gives
the same derivative gain also at integer exponents.
\end{remark}

Suppose $f\in \EEnd^\alpha(X)$ for some $\alpha\ge 1$, and let $\beta\in [0, \alpha-1]$. Analogously to the scalar duality in \eqref{def CL0}
on function spaces,
we define the transfer operator of $f$ from 
$C^\beta$ sections of $f^*TX$ to those of $TX$, that is, 
a bounded linear operator
$\ds \CL: \Gamma^\beta(f^*TX)\to \Gamma^\beta(TX)$
characterized by
\begin{equation}\label{eq: duality1}
\int_X
\lan(\CL v)(x),w(x)\rangle_x\,dm(x)
=
\int_X
\langle v(x),w(f(x))\rangle_{f(x)}\,dm(x)
\end{equation}
for every $v\in\Gamma^\beta(f^*TX)$ and
$w\in\Gamma^0(TX)$. 
We emphasize that the above pairings 
are evaluated
in the corresponding tangent spaces:
the pairing on the left is taken in $T_xX$ since
$(\CL v)(x),w(x)\in T_x X$, 
while the pairing on the right
is taken in $T_{f(x)}X$ since $v(x), w(f(x))\in T_{f(x)}X$.

We also need the transfer operator acting on covariant tensors. More precisely,
for any   $j\in \NN$, let
$\ds (T^*X)^{\otimes j}$
denote the bundle of covariant $j$-tensors on $X$.
We define the transfer operator on covariant $j$-tensor fields
along $f$, that is, a bounded linear operator
$\ds
\CL:\Gamma^\beta(f^*((T^*X)^{\otimes j}))\longrightarrow\Gamma^\beta((T^*X)^{\otimes j})
$
characterized by the following tensor duality:
\beq
\label{eq:tensor-duality}
\int_X
\left\langle(\CL A)(x), B(x)\right\rangle_x^{\mathrm{ct}}
\,dm(x)
=
\int_X
\left\langle A(x), B(f(x))\right\rangle_{f(x)}^{\mathrm{ct}}
\,dm(x)
\eeq
for every $A\in\Gamma^\beta(f^*((T^*X)^{\otimes j}))$ and
$B\in\Gamma^0((T^*X)^{\otimes j})$.
Here $\langle\cdot,\cdot\rangle_x^{\mathrm{ct}}$
denotes the inner product induced by the Riemannian metric
on the space $(T_x^*X)^{\otimes j}$ of covariant
$j$-tensors at $x$. 
 
\medskip 
 
Now we turn back to the setting of Theorem~\ref{thm: response eend}.
The following lemma shows that 
$\CL^{(k)}=\frac{1}{k!}\left. \p_t^k\CL_t\right|_{t=0}$
can be split into two parts.

\medskip

\begin{lemma}
\label{lem: explicit CLk}
 
For any $1\le k\le i\le r$  and any $\varphi\in \CB_i$, 
we have 
\beq
\label{explicit CLk}
\CL^{(k)}\varphi = -\dive\left( \CL\left(\varphi v_k\right)\right) + \CL^{(k)}_\#\varphi,
\eeq
where
$\CL^{(k)}_\#\varphi$ 
is determined by $\varphi$, 
$f$, $v_1, \dots, v_{k-1}$ but 
is independent of $v_k$.
Moreover, 
$\CL^{(k)}_\#\varphi$ has vanishing integral, i.e., $\ds \int \CL^{(k)}_\#\varphi dm=0$.
In particular, 
$\CL^{(1)}_\#\varphi=0$ 
and thus
$\ds
\CL^{(1)}\varphi = -\dive\left( \CL\left(\varphi v_1\right)\right).
$
\end{lemma}

\begin{proof}
Using the Riemannian exponential map, we write 
$\ds f_t(x)=\exp_{f(x)}\left( v(t, x)\right)$. 
For any $1\le k\le r$, 
the definition of $v_k$ given in \eqref{def vk} gives the following expansion:
\beq
\label{eq vtx k expansion}
v(t, x)=
t v_1(x)+t^2v_2(x)+\dots+t^k v_k(x)+o(|t|^k)\in T_{f(x)}X, \
\ \text{as} \ t\to 0.
\eeq
In particular, $\ds v(t, x)=t v_1(x) + o(|t|)$. 
For any $\psi\in C^k(X)$, the covariant Taylor expansion
yields
\beqn
\begin{split}
\psi(f_t(x))
&=\psi(f(x))+
\sum_{j=1}^k\frac{1}{j!}
\left(\nabla^j\psi \right)_{f(x)}
  \left(v(t,x),\ldots,v(t,x)\right) +o(|t|^k) \\
 &=\psi(f(x))+
\sum_{j=1}^k
 \left\langle \frac{1}{j!}
\bigl(v(t, x)^\flat\bigr)^{\otimes j},
(\nabla^j\psi)(f(x))
\right\rangle_{f(x)}^{\rct}
+o(|t|^k),
\quad\text{as }t\to 0,
\end{split}
\eeqn
where 
\begin{itemize}
\item
$\ds
v(t, x)^\flat\in T^*_{f(x)}X
$
is the metric-dual covector of $v(t, x)$, that is, 
$\ds v(t, x)^\flat (w)=\langle v(t, x), w\rangle_{f(x)}$
for any $w\in T_{f(x)}X$.
Hence
$\ds
\frac{1}{j!}
\bigl(v(t, x)^\flat\bigr)^{\otimes j}
$
is a covariant $j$-tensor at $f(x)$;
\item
$\nabla^j\psi$ denotes the $j$-th covariant derivative of $\psi$
with respect to the Levi-Civita connection,
which is a covariant $j$-tensor field.
Hence 
$\ds
(\nabla^j\psi)(f(x))
$
is also a covariant $j$-tensor at $f(x)$;
\item
The inner product $\langle\cdot,\cdot\rangle^\rct_{f(x)}$ 
is induced by the Riemannian metric  
on the space
$(T_{f(x)}^*X)^{\otimes j}$ of covariant $j$-tensors at $f(x)$.
\end{itemize}
Substituting the expansion of $v(t,x)$  in~\eqref{eq vtx k expansion}
and then collecting equal powers of $t$,
we rewrite 
\beqn
\frac{1}{j!}\bigl(v(t, x)^\flat\bigr)^{\otimes j}
=
\sum_{\ell=j}^{k} t^\ell A_{\ell, j}(x)+o(|t|^k),
\ \ \text{with} \
A_{\ell, j}(x):=\frac{1}{j!}
\sum_{\substack{k_1,\ldots, k_j\ge1\\
                k_1+\cdots+k_j=\ell}}
v_{k_1}(x)^\flat\otimes\cdots\otimes
v_{k_j}(x)^\flat.
\eeqn
Note that 
$A_{\ell, j}(x)$ is independent of $v_k(x)$
 for any $(\ell, j)\ne (k, 1)$.
Also,
$\ds
A_{k, 1}(x)= v_k(x)^\flat
$
and thus 
\beqn
\left\langle A_{k, 1}(x),
(\nabla^1\psi)(f(x))
\right\rangle_{f(x)}^\rct
=\left\langle v_k(x), \nabla \psi (f(x)) \right\rangle_{f(x)},
\eeqn
where 
$\nabla \psi$ in the second bracket
denotes the gradient vector field of $\psi$.
It follows that
\beq
\label{psi t Taylor}
\psi\left( f_t (x)\right) 
=\psi\left( f(x)\right) + 
t^k  \lan v_k(x), \nabla \psi( f(x))\ran_{f(x)}
+ \Psi_k (t, x)   +  o(|t|^k),
\ \ \text{as} \ t\to 0,
\eeq
where
\beq
\label{eq Psi k}
\Psi_k(t, x) = \sum_{\substack{1\le j\le \ell \le k \\ (\ell, j)\ne (k, 1)}}
 t^\ell 
 \left\langle  A_{\ell, j}(x),
(\nabla^j\psi)(f(x))
\right\rangle_{f(x)}^{\rct}
\eeq
is a polynomial in $t$ of degree  at most $k$,  with zero constant term.
We stress that once $f$, $\psi$ and the Riemannian metric are fixed,
$\Psi_k(t, x)$ depends only on $v_1, \dots, v_{k-1}$ but is independent of $v_k$.
In the particular case when $k=1$, we have $\Psi_1\equiv 0$ and thus
\beqn
\psi\left( f_t (x)\right) 
=\psi\left( f(x)\right) + 
t \lan v_1(x), \nabla \psi(f(x))\ran_{ f(x)} + o(|t|), \ \text{as} \ t\to 0.
\eeqn
We remark that all the above remainders $o(\cdot)$ are uniform in $x\in X$.

Now taking $\left.\frac1{k!}\partial_t^k\right|_{t=0}$ on both sides of the duality relation~\eqref{eq: duality},
and applying the expansion given by \eqref{psi t Taylor},
we obtain for any functions $\varphi\in\CB_i$ and $\psi\in C^k(X)$ that
\beqn
\int \CL^{(k)} \varphi(x) \psi(x) dm(x) 
=\int \varphi(x) \lan v_k(x), \nabla \psi(f(x))\ran_{f(x)}  dm(x)
+ \int  \varphi(x) \frac{\p_t^k \Psi_k\left(t, x\right)}{k!}\Big|_{t=0}   dm(x).
\eeqn
By \eqref{eq: duality1} and \eqref{dive Lib}, 
the first integral above  can be rewritten as 
\beqn
\begin{split}
\int   \lan \varphi(x) v_k(x), \nabla \psi(f(x))\ran_{f(x)}   dm(x)
&=\int   \lan \CL\left(\varphi v_k\right)(x), \nabla \psi(x) \ran_x  \, dm(x)  \\
&= -\int \psi(x) \, \dive\left(  \CL\left(\varphi v_k\right)\right)(x) \,  dm(x).
\end{split}
\eeqn
For the second integral, 
recall that $\Psi_k\left(t, x\right)$ is given by~\eqref{eq Psi k},
then  
\beqn
\left.\frac{\p_t^k \Psi_k\left(t, x\right)}{k!}\right|_{t=0}
=
\sum_{j=2}^{k}
 \left\langle  A_{k, j}(x),
(\nabla^j\psi)(f(x))
\right\rangle_{f(x)}^{\rct},
\eeqn
In particular, we have
$\ds
\left. \p_t \Psi_1\left(t, x\right)\right|_{t=0}\equiv 0
$
since $\Psi_1\equiv 0$.
By the tensor duality~\eqref{eq:tensor-duality} and repeated integration by parts, 
we rewrite the second integral as 
$\ds \int \CL^{(k)}_\#\varphi \psi dm$, where 
\beqn
\CL_\#^{(k)}\varphi:=
\sum_{j=2}^k
(\nabla^j)^*
\left[\CL(\varphi A_{k, j})\right].
\eeqn
Here $(\nabla^j)^*$ denotes the formal adjoint of $\nabla^j$.
It is easy to verify that $\CL_\#^{(k)}$
is a bounded linear operator from $\CB_i$ to $\CB_{i-k}$ for any $1\le k\le i\le r$.
Once the unperturbed map $f$ and the Riemannian metric are fixed,
$\CL_\#^{(k)}$ depends only on the first $(k-1)$ infinitesimal variations
$v_1,\ldots, v_{k-1}$ but it is independent
of $v_k$.
In particular, $\CL^{(1)}_\#\equiv 0$.

\medskip

To sum up, we obtain for any two  functions $\varphi\in\CB_i$ and $\psi\in C^k(X)$ that
\beq
\label{CLk explicit integral}
\int \CL^{(k)} \varphi \psi dm
=
- \int   \dive\left(  \CL\left(\varphi v_k\right)\right)  \psi dm
+ \int  \CL^{(k)}_\#\varphi \psi dm.
\eeq
Recall that $r\ge 1$ and $\al>r+2$,
then
$\CB_i:=C^{\alpha-r-1+i}(X) \subseteq C^k(X)$ for any $0\le k\le i\le r$.
The identity~\eqref{CLk explicit integral} holds for any $\varphi\in\CB_i$
and any $\psi\in C^k(X)$,
from which \eqref{explicit CLk} immediately follows.
Moreover, 
$\CL^{(k)}_\#\varphi$ has vanishing integral, since 
by~\eqref{CLj vanishing} and~\eqref{dive Lib0}, we have 
\beqn
 \int \CL^{(k)}_\#\varphi dm=
\int \CL^{(k)}\varphi dm + \int \dive\left( \CL\left(\varphi v_k\right)\right) dm=0.
\eeqn
The proof of Lemma~\ref{lem: explicit CLk} is complete.
\end{proof}

\subsection{Proof of Theorem~\ref{thm: response eend}}

We are now ready to prove Theorem~\ref{thm: response eend}, concerning the response formulae of any order.

\begin{proof}[Proof of Theorem~\ref{thm: response eend}] 
It is well known (see e.g.~\cite{MR1793194, MR2129258})  that
$\CP_t$
has the following representation:
\beqn
\CP_t\varphi =\rho_t \int \varphi dm,  
\ \ \text{for any} \ \
 \varphi\in \CB_i.
\eeqn
It follows  that 
$\rho_t=\CP_t 1 $,
where $1$ stands for the constant one function. 
Thus by Lemma~\ref{lem: der CP},   
\beqn
\rho_t - \rho = \left( \CP_t  -\CP\right) 1
=\sum_{j=1}^k t^j \CP^{(j)}1 + o\left( t^k\right),
\eeqn
where we take $i=r$ and the remainder is measured in
$C^{\alpha-k-2}(X)$.
The same expansions hold locally uniformly after recentering
at any nearby parameter. Consequently, the normalized $r$-th
forward differences converge locally uniformly in
$C^{\alpha-r-2}(X)$, identifying a continuous $r$-th
distributional derivative. Hence $t\mapsto\rho_t$ is a
$C^r$ curve in this space.
Thus
\beqn
\rho^{(k)}:= \left.\frac{1}{k!}\p_t^k \rho_t \right|_{t=0} = \CP^{(k)}1.
\eeqn
By Lemma~\ref{lem: CPkP} and 
the fact that $\CP1 =\rho$,  we further get
\beq
\label{rhok2}
\begin{split}
\rho^{(k)}
& =\sum_{(i_1 i_2\dots i_\ell) \in \Par(k)}  \CR \CL^{(i_1)}\CR \CL^{(i_2)} \dots \CR \CL^{(i_\ell)}  \rho \\
&=\CR\CL^{(k)}\rho
+\sum_{\substack{ (i_1 i_2\dots i_\ell) \in \Par(k): \\ \text{with} \  \ell\ge 2}} 
\CR \CL^{(i_1)}\CR \CL^{(i_2)} \dots \CR \CL^{(i_\ell)}\rho,
\end{split}
\eeq
in which $\CR$ is given by \eqref{def CR}.
By Lemma~\ref{lem: explicit CLk},
$\CL^{(j)}$ is independent of $v_k$ if $j<k$, 
and thus the summation in the second line of \eqref{rhok2} is independent of $v_k$.
Moreover, 
it follows from \eqref{explicit CLk} that the first term in the second line of \eqref{rhok2}
can be split into two parts, that is,
\beqn
\CR\CL^{(k)}\rho = - \CR\dive\left( \CL\left(\rho v_k\right)\right) + \CR \CL^{(k)}_\#\rho,
\eeqn
where the second term $\CR \CL^{(k)}_\#\rho$ is also independent of $v_k$.
Hence we obtain \eqref{explicit rhok} by setting  
\beq
\label{def rhokk}
\rho^{(k)}_\#= \CL^{(k)}_\#\rho + 
\sum_{\substack{ (i_1 i_2\dots i_\ell) \in \Par(k): \\ \text{with} \  \ell\ge 2}} 
\CL^{(i_1)}\CR \CL^{(i_2)} \dots \CR \CL^{(i_\ell)}\rho,
\eeq
which is determined by $f, v_1, \dots, v_{k-1}$ but is independent of $v_k$.
Due to \eqref{CLj vanishing} and Lemma~\ref{lem: explicit CLk},
all  terms in \eqref{def rhokk}
have vanishing integrals
and thus so does $\rho^{(k)}_\#$. 
In the particular case when $k=1$, we have
\beq
\label{def rho1'}
\rho^{(1)}=\CR\CL^{(1)}\rho =- \CR\dive\left( \CL\left(\rho v_1\right)\right). 
\eeq
The proof of Theorem~\ref{thm: response eend} is now complete.
\end{proof}

\section{Solution to Request Problem} 
\label{sec: proof2}

\subsection{\texorpdfstring{Auxiliary Operator $\CDD$}{Auxiliary Operator D}}

We first recall the assumptions of Theorem~\ref{thm: request}
and establish some basic properties.
More precisely, let $r\ge 1$ and $\alpha>r+2$, and let
$\left\{f; \wrho^{(1)}, \dots, \wrho^{(r)} \right\}$ be the prescribed data satisfying 
Assumptions \textbf{(R1)} and \textbf{(R2)} of Theorem~\ref{thm: request}.
The following statements hold:
\begin{enumerate}[(1)]
\item Let $\CL$ and $\rho$ denote the corresponding transfer operator and  ACIP density
of the unperturbed endomorphism $f$ respectively. 
By \eqref{def CL0} and that $\CL \rho =\rho$,
for any $\varphi\in L^1(m)$ and $\psi\in L^\infty(m)$, 
\beqn
\int \CL\left(\rho \, \varphi\circ f \right) \psi dm
=
\int  \rho \, \varphi\circ f \, \psi\circ f dm
= 
\int \CL \rho \, \varphi \, \psi dm
=
\int \rho \, \varphi  \, \psi dm,
\eeqn
which implies 
the following equality:
\beq
\label{def normal transfer}
\CL\left(\rho \, \varphi\circ f \right)=\rho \, \varphi,
\ \ \text{for any} \ \varphi\in L^1(m).
\eeq
Note that \eqref{def normal transfer} also holds when $\varphi$ is
a vector field in $\Gamma^\beta(TX)$ for any $\beta\in [0, \alpha+r-1]$.
\item 
We assume that $\wrho^{(k)}\in C^{\alpha+r-1-k}(X)$ for any $1\le k\le r$,
such that $\ds \int \wrho^{(k)} dm=0$. Then  
by taking $\varphi=\wrho^{(k)}$ and $\psi\equiv 1$ in \eqref{def CL0}, we get
$\ds \int \CL\wrho^{(k)} dm =0$, and thus
\beq
\label{vanishing1}
\int (I-\CL) \wrho^{(k)} dm=
\int \wrho^{(k)} dm - \int \CL\wrho^{(k)} dm =0,
\eeq
that is, 
$(I-\CL) \wrho^{(k)}$ also has vanishing integral.
\end{enumerate}

To solve the request problem, 
we define an operator $\CDD$ acting on a section $v$ of $f^*TX$ by
\beq
\label{def CDD}
\CDD v:=-\dive\left( \CL \left( \rho v \right)\right),
\eeq
where $\dive(\cdot)$ is the divergence operator given in \eqref{def divergence}.

\begin{lemma}
\label{lem: CDD}
The operator $\CDD$ is a bounded linear operator from $\Gamma^{\beta}(f^*TX)$
to $C^{\beta-1}(X)$ for any $\beta\in [1, \alpha+r-1]$,
such that the following statements hold:
\begin{enumerate}[(1)]
\item $\ker(\CDD)$ is an infinite dimensional subspace of $\Gamma^{\beta}(f^*TX)$;
\item Assume additionally that $\beta>1$ is noninteger, or that $\dim X=1$. For any $\varphi \in C^{\beta-1}(X)$ such that $\ds \int \varphi dm=0$,
the equation $\CDD v =\varphi$ has general solutions of the form 
$v=\xi+\hv$, where 
$\xi\in \ker(\CDD)$ and
$\hv\in \Gamma^{\beta}(f^*TX)$ is a particular solution.
\end{enumerate}
\end{lemma}

\begin{proof}
The boundedness of $\CDD$ from $\Gamma^{\beta}(f^*TX)$
to $C^{\beta-1}(X)$ is due to the facts that 
$\rho\in C^\beta(X)$, as well as that
$\CL: \Gamma^{\beta}(f^*TX) \to \Gamma^{\beta}(TX)$
and 
$\dive: \Gamma^{\beta}(TX)\to C^{\beta-1}(X)$
are both bounded linear operators. 

(1) We claim that 
$\ker\left(\CL \right)$ is an infinite dimensional subspace of $\Gamma^{\beta}(f^*TX)$.
Indeed, the transfer operator $\CL$ has an expression given by 
\eqref{def transfer2}, and thus for any $C^\beta$ section $w_1\in \Gamma^{\beta}(f^*TX)$  that is 
compactly supported in one inverse branch, it is easy to extend $w_1$  to some
$w\in\Gamma^{\beta}(f^*TX)$ such that $\CL(w)=0$.
Since the space of $C^\beta$ sections compactly supported in a nonempty open set is infinite dimensional,
so is $\ker\left(\CL \right)$ (see Remark~\ref{rem: solution infinite} for details).
It follows from~\eqref{def CDD} that
$\ker(\CDD)$ contains $\rho^{-1} \ker\left(\CL \right)$,
and thus $\ker(\CDD)$ is also an infinite dimensional subspace of $\Gamma^{\beta}(f^*TX)$.

(2)  
Now given $\varphi \in C^{\beta-1}(X)$ such that $\ds \int \varphi dm=0$,
by \eqref{div eq}, there exists 
$w\in \Gamma^{\beta}(TX)$ such that
$\ds\dive(w)=\varphi$.
It follows from~\eqref{def normal transfer}  that 
$\ds
\hv:=-\left(\rho^{-1} w\right)\circ f
$
is a particular solution to the equation $\CDD v =\varphi$,
since
\beqn
\CDD \hv 
=\dive\left(\CL \left( \rho \left(\rho^{-1} w\right)\circ f\right) \right)
=\dive\left( \rho \cdot \rho^{-1} w\right) =\dive(w)=\varphi.
\eeqn
Thus the general solution of $\CDD v =\varphi$ is given 
by $v=\xi+\hv$ for any $\xi\in \ker(\CDD)$.
Here $\hv(x)=-\rho(f(x))^{-1}w(f(x))\in T_{f(x)}X$, so $\hv$ is indeed a section of $f^*TX$.
\end{proof}

\begin{remark}
\label{rem: solution infinite}
Since $f\in \EEnd^{\alpha+r}(X)$ is a covering with at least two sheets,
we choose an evenly covered open set $U$ and two inverse branches
$g_1,g_2:U\to X$ with disjoint images $U_1,U_2$.
Write $J_f=|\det Df|$.
For any $w_1\in \Gamma^\beta(f^*TX)$ that is compactly supported
in $U_1$, set $w=w_1$ on $U_1$ and define
\[
w(g_2(y))
=
-\frac{J_f(g_2(y))}{J_f(g_1(y))}
w_1(g_1(y)),
\quad \text{for any} \ y\in U.
\]
Set $w=0$ outside $U_1\cup U_2$.
Both $w_1(g_1(y))$ and $w(g_2(y))$ belong to $T_yX$,
so this definition is intrinsic.
The regularity of $J_f$ and the compact support of $w_1$
ensure that $w\in\Gamma^\beta(f^*TX)$.
For any $y\in U$, the only potentially nonzero branch contributions
cancel, that is,
\[
(\CL w)(y)
=
\frac{w_1(g_1(y))}{J_f(g_1(y))}
+
\frac{w(g_2(y))}{J_f(g_2(y))}
=0.
\]
For $y\notin U$, all branch contributions vanish.
Thus $\CL w=0$.
The construction is linear and injective in $w_1$, since
$w|_{U_1}=w_1$.
A local trivialization and linearly independent smooth bump
functions supported in $U_1$ therefore yield infinitely many
linearly independent elements of $\ker\CL$.
\end{remark}

\subsection{Proof of Theorem~\ref{thm: request}}

We are now ready to prove Theorem~\ref{thm: request},
concerning the existence of solutions to the request problem - Problem~\ref{qu: request}.

\begin{proof}[Proof of Theorem~\ref{thm: request}]

Recall that the prescribed data
$\left\{f; \wrho^{(1)}, \dots, \wrho^{(r)} \right\}$ 
satisfies Assumptions \textbf{(R1)} and \textbf{(R2)},
and denote by
 $\CV\left(f; \wrho^{(1)}, \dots, \wrho^{(r)}\right)$ 
the solution set to the request problem.
It follows from Theorem~\ref{thm: response eend} 
that any
$\ds (v_1, \dots, v_r)\in \CV\left(f;  \wrho^{(1)}, \dots, \wrho^{(r)}\right)$
is a solution to the following equation system:

\beq
\label{solve v}
\begin{cases}
\wrho^{(1)}&= \ \CR\, \left[ -\dive\left( \CL\left(\rho v_1\right)\right)\right], \\
\wrho^{(2)}&= \ \CR\, \left[- \dive\left( \CL\left(\rho v_2\right)\right) + \rho^{(2)}_\#(v_1)\right], \\
 & \vdots \\
\wrho^{(k)}&= \ \CR\, \left[ -\dive\left( \CL\left(\rho v_k\right) \right) + \rho^{(k)}_\#(v_1, \dots, v_{k-1})\right], \\
& \vdots\\
\wrho^{(r)}&= \ \CR\, \left[ - \dive \left(\CL\left(\rho v_r\right) \right) + \rho^{(r)}_\#(v_1, \dots, v_{r-1})\right].
\end{cases}
\eeq
Here we write  
$\ds
\rho^{(k)}_\#=\rho^{(k)}_\#(v_1, \dots, v_{k-1})
$
to emphasize the dependence of $\ds \rho^{(k)}_\#$ on $\{v_1, \dots, v_{k-1}\}$.

\medskip

Recall that $\CR$ is given by \eqref{def CR}, and it is easy to verify that 
$(I-\CL)\CR=I-\CP$.
We also notice that all the brackets $[\dots]$ on the right hand side of \eqref{solve v}
have vanishing integrals, then by \eqref{def CP0},
the action of $\CP$ on all brackets $[\dots]$ vanishes.
Hence we get
\beqn
(I-\CL)\CR [\dots] = [\dots]
\eeqn
for all the brackets $[\dots]$ in \eqref{solve v}.
We also recall the definition of the operator $\CDD$ given by \eqref{def CDD}.
Then~\eqref{solve v} is equivalent to the following equation system:

\beq
\label{solve v1}
\begin{cases}
(I-\CL)\wrho^{(1)}&= \ \CDD v_1, \\
(I-\CL)\wrho^{(2)}&=  \ \CDD v_2+ \rho^{(2)}_\#(v_1), \\
 & \vdots \\
(I-\CL)\wrho^{(k)}&=  \ \CDD v_k + \rho^{(k)}_\#(v_1, \dots, v_{k-1}), \\
& \vdots\\
(I-\CL)\wrho^{(r)}&= \  \CDD v_r + \rho^{(r)}_\#(v_1, \dots, v_{r-1}).
\end{cases}
\eeq

\medskip

We can then solve the equation system~\eqref{solve v1} step by step in the class
$v_k\in\Gamma^{\alpha+r-k}(f^*TX)$, $1\le k\le r$, as follows.

\medskip 

\noindent
$\bullet$\,
Since $\wrho^{(1)}\in C^{\alpha+r-2}(X)$ satisfies 
$\ds \int (I-\CL)\wrho^{(1)}=0$ (see \eqref{vanishing1}), 
by Lemma~\ref{lem: CDD}, 
the first equation 
\beqn
\CDD v_1 =(I-\CL)\wrho^{(1)}
\eeqn
has general solutions of the form $v_1=\xi_1+ \hv_1 $,
where $\ds \xi_1\in \ker\left( \CDD|_{\Gamma^{\alpha+r-1}(f^{*}TX)}\right)$ and
$\hv_1\in \Gamma^{\alpha+r-1}(f^{*}TX)$.
Here we denote the particular solution by $\hv_1=\hv_1\left(\wrho^{(1)}\right)$
to emphasize its dependence on $\wrho^{(1)}$.

\medskip

\noindent
$\bullet$\,
Suppose now the first $(k-1)$ equations in~\eqref{solve v1} have already been solved, where $2\le k\le r$.
That is, for any $1\le j\le k-1$, the $j$-th equation of~\eqref{solve v1} has general solutions of 
the form $v_j=\xi_j+\hv_j $, where $\ds \xi_j\in \ker\left( \CDD|_{\Gamma^{\alpha+r-j}(f^{*}TX)}\right)$ and
$\ds
\hv_j=\hv_j\left(\wrho^{(1)}, \dots, \wrho^{(j)}, \xi_{1}, \dots, \xi_{j-1}\right)\in \Gamma^{\alpha+r-j}(f^{*}TX).
$
Then the $k$-th equation of~\eqref{solve v1} becomes an equation with only one unknown variable $v_k$, that is,
\beq
\label{def k-th eq}
\CDD v_k=(I-\CL)\wrho^{(k)} - \rho^{(k)}_\#(v_1, \dots, v_{k-1}),
\eeq
of which the right hand side has vanishing integral due to \eqref{vanishing1} and 
Item (3) of Theorem~\ref{thm: response eend}.
Grouping the ordered partitions in \eqref{def rhokk} by their
first entry and using the previously solved request equations,
we obtain
\[
\rho^{(k)}_\#=\CL^{(k)}_\#\rho+
\sum_{j=1}^{k-1}\CL^{(j)}\wrho^{(k-j)}.
\]
It follows that
$\rho^{(k)}_\#\in C^{\alpha+r-1-k}(X)$.
Thus by Lemma~\ref{lem: CDD} again, 
Equation~\eqref{def k-th eq}  has general solutions of the form $v_k=\xi_k +\hv_k$,
where $\ds \xi_k\in \ker\left( \CDD|_{\Gamma^{\alpha+r-k}(f^{*}TX)}\right)$
and
\beqn
\hv_k=\hv_k\left(\wrho^{(1)}, \dots, \wrho^{(k)}, \xi_{1}, \dots, \xi_{k-1}\right)\in \Gamma^{\alpha+r-k}(f^{*}TX).
\eeqn

To sum up, the solutions in
$\ds \CV\left(f; \wrho^{(1)}, \dots, \wrho^{(r)}\right)$
satisfying
\beq
\label{solution consistent}
(v_1, \dots, v_r) \in 
\prod_{k=1}^r  \Gamma^{\alpha+r-k}(f^{*}TX)  \subseteq 
\left[{\Gamma^\alpha(f^{*}TX)}\right]^r
\eeq
are precisely the tuples of the form
\beq
\label{explicit request sol}
\begin{cases}
v_1&= \ \xi_1 + \hv_1\left( \wrho^{(1)} \right), \\
v_2&= \ \xi_2 + \hv_2\left(\wrho^{(1)},   \wrho^{(2)}, \xi_{1}\right), \\
 & \vdots \\
v_k&= \ \xi_k + \hv_k\left(\wrho^{(1)}, \dots, \wrho^{(k)}, \xi_{1}, \dots, \xi_{k-1}\right), \\
& \vdots\\
v_r &= \ \xi_r + \hv_r\left(\wrho^{(1)}, \dots, \wrho^{(r)}, \xi_{1}, \dots, \xi_{r-1}\right).
\end{cases}
\eeq
Every constructed tuple  
 is realized by the $C^r$ parameter family
\(\ds
f_t(x)=\exp_{f(x)}
\left(\sum_{k=1}^r t^k v_k(x)\right),
\)
which consists of expanding endomorphisms for sufficiently small $|t|$.
Moreover, fixing $v_1,\ldots,v_{r-1}$ and varying $v_r$ by elements
of the infinite-dimensional space
$\ker(\CDD|_{\Gamma^\alpha(f^*TX)})$
preserves all request equations, and thus yields an infinite-dimensional
affine subspace of the solution set.
The proof of Theorem~\ref{thm: request} is complete.
\end{proof}

\subsection{Proof of Corollary~\ref{cor: request opt}}

Corollary~\ref{cor: request opt} shows that there exist minimizers of the relaxed request problem.

\begin{proof}[Proof of Corollary~\ref{cor: request opt}]
Note that the space $\left[{\Gamma^{\alpha}(f^*TX)}\right]^r$
is compactly embedded in $\left[\Gamma^{\beta}(f^*TX)\right]^r$
such that $\|\bv\|_\beta\le  \|\bv\|_\alpha$
 for any 
$\bv \in \left[{\Gamma^{\alpha}(f^*TX)}\right]^r$.
Therefore, for any $K>0$, the set
\beqn
\CV_K:=\left\{\bv\in  \CV: \ 
\|\bv\|_\alpha\le K\right\}
\eeqn
is a bounded subset of $\left[\Gamma^{\alpha}(f^*TX)\right]^r$,
and thus $\overline{\CV_K}^{\,\beta}$ is a compact subset of $\left[\Gamma^{\beta}(f^*TX)\right]^r$.
By Condition (1), we may choose 
a sufficiently large $K>0$ such that 
\beqn
\inf\left\{\CJ(\bv): \ \bv\in  \CV\backslash \CV_K\right\}>
\inf\left\{\CJ(\bv): \ \bv\in  \CV \right\}.
\eeqn
Hence $\inf_{\CV_K}\CJ=\inf_{\CV}\CJ$.
By its definition, $\overline{\CJ}$ is lower semicontinuous
on $\overline{\CV}^{\,\beta}$, and Condition (2) gives
$\overline{\CJ}=\CJ$ on $\CV$.
Thus, there must exist  
 $\overline{\bv}\in\overline{\CV_K}^{\,\beta}\subseteq\overline{\CV}^{\,\beta}$ such that
\beqn
\inf_{\CV}\CJ\le\overline{\CJ}(\overline{\bv})
=\min_{\bv\in\overline{\CV_K}^{\,\beta}}\overline{\CJ}(\bv)
\le\inf_{\bv\in\CV_K}\CJ(\bv)=\inf_{\CV}\CJ.
\eeqn
Thus $\overline{\bv}$ also minimizes
$\overline{\CJ}$ on $\overline{\CV}^{\,\beta}$.
The proof of Corollary~\ref{cor: request opt} is now complete.
\end{proof}


\section{Applications to Circle Expanding Endomorphisms}

In \cite{MR3670005}, Galatolo and Pollicott  presented 
a  solution to the linear request problem 
in the particular case of circle expanding endomorphisms.
In this section, we  work with the same category and
apply our main theorems to obtain 
explicit  solutions to the quadratic request problem.

In Subsection~\ref{subsec: example general}, we present  explicit formulae for constructing perturbations with prescribed first and second order responses. The second order equation depends quadratically on the first order perturbation, which makes the request constraints nonlinear. Consequently, changing the first order perturbation may require a corresponding second order correction, even if the first order response remains unchanged.

In Subsection~\ref{subsec: opt quad doubling}, we consider the doubling map.
A perturbation can have a zero first order response but a nonzero second order response. 
Thus, the linear response alone does not capture its statistical effect.
Finally, Example~\ref{example: quadratic non-uniq} shows that
the relaxed minimizer need not be unique.

\subsection{General Solutions to Quadratic Request on the Circle}
\label{subsec: example general}

Let $r=2$ and $\al=5$.
The perturbation families below are $C^2$ curves in $\EEnd^5(\TT)$;
their Taylor remainders are taken in $C^5(\TT)$, and the density
remainders are taken in $C^1(\TT)$.
We consider the class $\EEnd^5(\TT)$ of circle expanding endomorphisms,
where $\TT\cong \RR/\ZZ$ is the unit circle.
Functions and vector fields on $\TT$ are represented by one-periodic functions on $\RR$. A circle map of degree $d$ has a lift $F$ satisfying $F(x+1)=F(x)+d$; its perturbations below are written using such lifts.
The quadratic request problem in the  setting of 
circle expanding endomorphisms 
is then restated as follows.

\begin{problem}[Quadratic Request for $\EEnd^5(\TT)$]  
\label{Pro: Q}
Given    $f\in \EEnd^5(\TT)$, $\wrho^{(1)}\in C^3(\TT)$ and $\wrho^{(2)}\in C^2(\TT)$,
we seek $(v_1, v_2)\in C^5(\TT)\times  C^5(\TT)$ 
such that the following holds:
if  a family of circle expanding endomorphisms
$\{f_t\}\subseteq \EEnd^5(\TT)$ has an expansion of the form
\beqn
f_t= f+t v_1 + t^2 v_2 + o(t^2),
\ \ \text{as} \ t\to 0,
\eeqn
then the corresponding family of ACIP densities $\{\rho_t\}$ has an expansion of the form
\beqn
\rho_t=\rho+t\wrho^{(1)} + t^2\wrho^{(2)} + o(t^2),
\ \ \text{as} \ t\to 0.
\eeqn

The set of all such pairs $(v_1, v_2)$ is called 
the solution set of the quadratic request problem in the category  $\EEnd^5(\TT)$,
denoted as
$\ds \CV\left(f; \wrho^{(1)},  \wrho^{(2)}\right)$.
\end{problem}

\medskip

As shown in Theorem~\ref{thm: request},
Assumptions \textbf{(R1)} and \textbf{(R2)} are imposed on
the prescribed data $\{f; \wrho^{(1)}, \wrho^{(2)}\}$,
that is,
\beq
\label{Assump R circle}
f\in \EEnd^7(\TT), \ \
\wrho^{(1)}\in C^5(\TT), \ \
\wrho^{(2)}\in C^4(\TT), 
\ \text{and} \ \
\int \wrho^{(1)}(x) \, dx = \int \wrho^{(2)}(x) \, dx=0.
\eeq
 

\begin{prop}
\label{prop: qu req}
Under Condition~\eqref{Assump R circle},
the solution set $\CV\left( f; \wrho^{(1)}, \wrho^{(2)} \right)$ 
to the quadratic request problem, i.e., Problem~\ref{Pro: Q},
contains all
the pairs $(v_1, v_2)$ of the form
\beq
\label{sol qu req}
\begin{cases}
v_1= \xi_1 + \left(\rho^{-1} w_1 \right)\circ f, \\
v_2= \xi_2 + \left(\rho^{-1} \left(w_2 + w_{2\#}  \right)\right)\circ f,
\end{cases}
\eeq
where 
\begin{enumerate}[(1)]
\item 
$\rho$ is the ACIP density function of $f$;
\item
$\ds
\xi_k\in \rho^{-1} \ker\left( \left. \CL\right|_{C^{7-k}(\TT)} \right)
$
for $k=1, 2$,
in which $\CL$ is the transfer operator of $f$;
\item
$w_k$ is an anti-derivative of $(\CL-I)\wrho^{(k)}$ for $k=1, 2$;
\item
$w_{2\#}$ is given by the following formula:
\beq
\label{def w22}
w_{2\#}:=\frac{1}{2}\left(\CL\left(\rho v_1^2\right)\right)'
- \CL\left(\wrho^{(1)} v_1\right).
\eeq
\end{enumerate}
\end{prop}

\medskip

\begin{proof}
It follows from \eqref{solve v1} that the
quadratic request problem boils down to solving 
the following equation system:
\beq
\label{solve v2}
\begin{cases}
(I-\CL)\wrho^{(1)}&= \  \CDD v_1, \\
(I-\CL)\wrho^{(2)}&= \  \CDD v_2+ \rho^{(2)}_\#,
\end{cases}
\eeq
where 
the operator $\CDD$ is defined by \eqref{def CDD},
and it can be rewritten as 
$
\CDD v:=-\left( \CL \left( \rho v \right)\right)',
$
since the divergence is exactly the derivative (with respect to the phase variable $x$)
for unit circle; furthermore,
the term $\rho^{(2)}_\#$ is given by 
\beqn
\rho^{(2)}_\# =\frac{1}{2}\left(\CL\left(\rho v_1^2\right)\right)''
-\left( \CL\left(\wrho^{(1)} v_1\right)\right)'.
\eeqn
Indeed, for any $\psi\in C^2(\TT)$, as $t\to 0$, we have  
\beqn
\begin{split}
\psi\circ f_t
&=\psi\circ \left(f + tv_1 + t^2v_2  \right) + o(t^2) \\
&= \psi\circ f + tv_1 \cdot \psi'\circ f  +
t^2\left(\frac12 v_1^2 \cdot  \psi'' \circ f   + v_2 \cdot  \psi'\circ f   \right)+ o(t^2).
\end{split}
\eeqn
Taking derivatives on both sides of the  duality relation~\eqref{eq: duality},\
we immediately get 
\beqn
\CL^{(1)}\varphi=-\left( \CL\left(\varphi v_1\right)\right)', \ \ 
\text{and} \ \ 
\CL^{(2)}\varphi=-\left( \CL\left(\varphi v_2\right)\right)' +\frac12 \left( \CL\left(\varphi v_1^2\right)\right)''.
\eeqn
By  Lemma~\ref{lem: explicit CLk}, we have
$\ds \CL_\#^{(2)}\varphi=\frac12 \left( \CL\left(\varphi v_1^2\right)\right)''$. 
Then by \eqref{def rhokk} and~\eqref{def rho1'}, we get
\beqn
\rho^{(2)}_\# =\CL^{(2)}_\#\rho + 
 \CL^{(1)}\CR \CL^{(1)}  \rho
 =\CL^{(2)}_\#\rho +\CL^{(1)} \wrho^{(1)}
 =\frac{1}{2}\left(\CL\left(\rho v_1^2\right)\right)''
-\left( \CL\left(\wrho^{(1)} v_1\right)\right)'.
\eeqn

Note that
the equation system \eqref{solve v2} is further converted into 
\beq
\label{solve v3}
\begin{cases}
\CL \left( \rho v_1 \right) = w_1, \\
\CL \left( \rho v_2 \right) =  w_2+w_{2\#},
\end{cases}
\eeq
where $w_k$ is an anti-derivative of $(\CL-I)\wrho^{(k)}$ for $k=1, 2$, and
$w_{2\#}$ is the anti-derivative of $\rho^{(2)}_\#$ given by \eqref{def w22}.
By \eqref{def normal transfer},
the equation system \eqref{solve v3} has a particular solution
\beqn
\begin{cases}
\hv_1=  \left(\rho^{-1} w_1 \right)\circ f, \\
\hv_2=  \left(\rho^{-1} \left(w_2 + w_{2\#}  \right)\right)\circ f.
\end{cases}
\eeqn
Therefore, all pairs given by \eqref{sol qu req}
with the stated regularity solve \eqref{solve v3}.
\end{proof}

\subsection{Optimal Quadratic Request for  Doubling Map}
\label{subsec: opt quad doubling}

We now consider a simple example  
which can be analyzed using the classical Fourier analysis on the unit circle.
To be precise, 
let $f: \TT\to \TT$ be the doubling map defined by $f(x)=2x \pmod 1$, and 
then the corresponding transfer operator is given by
\beqn
\label{def trans op doubling}
\CL\varphi(x)=\frac12\left( \varphi\left(\frac{x}{2}\right) + \varphi\left(\frac{x+1}{2}\right) \right),
 \ \ \text{for any} \ \varphi\in L^1(\TT, dx),
\eeqn
where $dx$ stands for the Lebesgue measure on $\TT$.
It is easy to see that 
the ACIP of  $f$ is exactly the Lebesgue measure, and thus
its ACIP density is $\rho\equiv 1$.

\medskip

Let $\{e_n\}_{n\in \ZZ}$ be
the canonical orthonormal Fourier basis on $\TT$ under the $L^2$ inner product,
where $e_n(x):=e^{i 2\pi n x}$. 
Consider two linear subspaces of $L^2(\TT, dx)$ given as follows:
\beqn
\CB_{\mathrm{odd}}:=\overline{\mathrm{Span}\left\{e_n: \ n\ \text{is  odd}\right\}}
\ \ \text{and} \ \ 
\CB_{\mathrm{even}}:=\overline{\mathrm{Span}\left\{e_n: \ n\ \text{is  even}\right\}},
\eeqn
that is, 
$\CB_{\mathrm{odd}}$ (or $\CB_{\mathrm{even}}$) consists
of $L^2$ functions with vanishing even (or odd) Fourier coefficients.

\smallskip

\begin{lemma} 
\label{lem: CB doubling}
Let $\CL$ be the transfer operator of doubling map $f$. 
The following statements hold:
\begin{enumerate}[(1)]
\item
Every function $\varphi\in L^2(\TT, dx)$ can be decomposed as 
$\varphi=\varphi_{\mathrm{odd}}+\varphi_{\mathrm{even}}\in 
\CB_{\mathrm{odd}}\oplus \CB_{\mathrm{even}}$, and
by the orthogonality, we have
$
\|\varphi\|_{L^2}^2=\|\varphi_{\mathrm{odd}}\|_{L^2}^2+\|\varphi_{\mathrm{even}}\|_{L^2}^2
$.
\item 
Let ``$\cdot$" denote the multiplication of two functions.
In the following inclusions, the factors are assumed to
belong to $L^\infty$ as well as to the indicated subspaces, so their
products belong to $L^2$. This applies to all later uses.
We have that
$\CB_{\mathrm{odd}}\cdot \CB_{\mathrm{odd}}\subseteq \CB_{\mathrm{even}}$,
$\CB_{\mathrm{even}}\cdot \CB_{\mathrm{even}}\subseteq \CB_{\mathrm{even}}$
and 
$\CB_{\mathrm{odd}}\cdot \CB_{\mathrm{even}}\subseteq \CB_{\mathrm{odd}}$.
\item
$\left(C^1(\TT)\cap\CB_{\mathrm{odd}}\right)'\subseteq \CB_{\mathrm{odd}}$
and
$\left(C^1(\TT)\cap \CB_{\mathrm{even}}\right)'\subseteq \CB_{\mathrm{even}}$,
where the prime represents the  derivative 
with respect to the phase variable $x$.
Also, $\left(\CL\varphi\right)'=\frac12 \CL\varphi'$ for any $\varphi\in C^1(\TT)$.
\item 
$
\CL e_n = 
\begin{cases}
e_{\frac{n}{2}},  &  \text{if}\ n \ \text{is even}, \\
0,   &  \text{if}\ n \ \text{is  odd}, \\
\end{cases}
$
and thus
$
\ker(\CL)=\CB_{\mathrm{odd}}.
$
\item
Since $e_n\circ f=e_{2n}$, we have
$\varphi\circ f\in \CB_{\mathrm{even}}$ for any $\varphi\in L^2(\TT, dx)$.
Also, $\left(\varphi\circ f\right)'=2\varphi'\circ f$
for any $\varphi\in C^1(\TT)$.
\item
For any $\varphi\in L^2(\TT, dx)$, we have
$\ds
\CL(\varphi\circ f)=\varphi
$
and
$\ds
\CL\varphi\circ f=\varphi_{\mathrm{even}}
$.
\end{enumerate}
\end{lemma}

\smallskip

The proof of this lemma
involves simple verifications, and is
left to the readers.

\medskip

Now we provide general solutions to the quadratic request problem
in the case when the unperturbed map $f$ is the doubling map.

\begin{prop}
\label{prop: qu req double}
Let $f$ be the doubling map of $\TT$.
Assume that $\wrho^{(1)}$ and $\wrho^{(2)}$ satisfy 
Condition~\eqref{Assump R circle}, and 
their Fourier series expansions are given by
\beqn
\wrho^{(k)}=\sum_{n\in \ZZ} \wrho_{k, n}e_n,
 \ \ \text{for} \ \ k=1, 2.
\eeqn
Then the solution set $\CV\left( f; \wrho^{(1)}, \wrho^{(2)} \right)$
to the quadratic request problem
contains all
the pairs $(v_1, v_2)$ of the form
\beqn
\begin{cases}
v_1= c_1 + \xi_1 +   \hv_1, \\
v_2= c_2 + \xi_2 +  \hv_2 + \CU(c_1, \xi_1),    
\end{cases}
\eeqn
where 
\begin{enumerate}[(1)]
\item free variables: 
$\ds c_k\in \RR$ 
and 
$\ds
\xi_k\in  C^{7-k}(\TT)\cap \CB_{\mathrm{odd}}
$
for $k=1, 2$;
\item
special solutions:
$\ds
\begin{cases}
\hv_1=  \ww_1\circ f, \\
\hv_2= (\ww_2-\ww_1\wrho^{(1)})\circ f,    
\end{cases}
$
where $\ww_k$ is given by
\beq
\label{def wwk}
\ww_k:=\sum_{n\in \ZZ\backslash\{0\}} \frac{\wrho_{k, 2n} - \wrho_{k, n}}{i 2\pi n} \, e_n,
\ \ \text{for}\ k=1, 2.
\eeq
\item 
corrected term in $v_2$: $\CU(c_1, \xi_1)$ is given by
\beq
\label{def CU}
\CU(c_1, \xi_1)= 
\frac12 \xi_1\xi_1'  - \wrho^{(1)}_{\mathrm{odd}} \xi_1 - c_1\wrho^{(1)} \circ f \in \CB_{\mathrm{even}}.
\eeq
\end{enumerate}
\end{prop}

\begin{proof} By Proposition~\ref{prop: qu req}, we only need to compute the  terms in \eqref{sol qu req}
for the case of doubling map. 
Note that $\rho\equiv 1$, and by Item (4) of
Lemma~\ref{lem: CB doubling}, we have
$\ds \xi_k\in \ker(\CL|_{C^{7-k}(\TT)})=C^{7-k}(\TT)\cap \CB_{\mathrm{odd}}$ for $k=1, 2$.
Moreover, the anti-derivatives  of 
\beqn
(\CL-I)\wrho^{(k)}=\CL \wrho^{(k)}- \wrho^{(k)}
=\sum_{n\in \ZZ } \left(\wrho_{k, 2n} - \wrho_{k, n}\right) e_n
\eeqn
are given by
$w_k=c_k + \ww_k$,
where $c_k\in \RR$ and
$\ww_k$ is given by \eqref{def wwk}.
It remains to calculate  $w_{2\#}\circ f$.
Recall that $w_{2\#}$ is given by \eqref{def w22}, then
by Item (3) (6) of
Lemma~\ref{lem: CB doubling}, we have
 \beqn
w_{2\#} \circ f
=\left(\frac{1}{2}\left(\CL\left( v_1^2\right)\right)' 
- \CL\left(\wrho^{(1)} v_1\right)    \right)\circ f
=\CL \left( \frac12 v_1v_1 ' 
-  \wrho^{(1)} v_1  \right)\circ f
=\left( \frac12 v_1v_1 ' 
-  \wrho^{(1)} v_1  \right)_{\mathrm{even}}.
\eeqn
Since 
$\ds 
v_1=\xi_1 + w_1 \circ f
$,
by Item (5) of
Lemma~\ref{lem: CB doubling},
we have 
$\ds 
v_1'=\xi_1' + 2w_1' \circ f
$,
and thus
\beqn
\begin{split}
& \frac12 v_1v_1 ' = \left(\frac12 \xi_1\xi_1 '  + \left(w_1 w_1'\right)\circ f\right)
+ \left( \xi_1 w_1'\circ f + \frac12 \xi_1' w_1\circ f\right)
\in \CB_{\mathrm{even}} \oplus \CB_{\mathrm{odd}};
\\
& \wrho^{(1)} v_1 = 
\left( \wrho^{(1)}_{\mathrm{odd}} \xi_1 +  \wrho^{(1)}_{\mathrm{even}}  w_1\circ f \right) 
+ \left( \wrho^{(1)}_{\mathrm{even}} \xi_1 +  \wrho^{(1)}_{\mathrm{odd}}  w_1\circ f \right) 
\in \CB_{\mathrm{even}} \oplus \CB_{\mathrm{odd}}.
\end{split}
\eeqn
Therefore, 
\beqn
\begin{split}
w_{2\#} \circ f 
& =\frac12 \xi_1\xi_1 '  + \left(w_1 w_1'\right)\circ f 
- \wrho^{(1)}_{\mathrm{odd}} \xi_1 -  \wrho^{(1)}_{\mathrm{even}}  w_1\circ f \\
& = \left( \frac12 \xi_1\xi_1 '  - \wrho^{(1)}_{\mathrm{odd}} \xi_1\right)
+w_1\circ f \left( w_1'\circ f - \wrho^{(1)}_{\mathrm{even}}\right) \\
& = \left( \frac12 \xi_1\xi_1 '  - \wrho^{(1)}_{\mathrm{odd}} \xi_1\right)
- (w_1 \wrho^{(1)})\circ f \\
&= \left( \frac12 \xi_1\xi_1 '  - \wrho^{(1)}_{\mathrm{odd}} \xi_1
- c_1\wrho^{(1)} \circ f \right)  - (\ww_1 \wrho^{(1)})\circ f \\
&=\CU(c_1, \xi_1) - (\ww_1 \wrho^{(1)})\circ f,
\end{split}
\eeqn
where $\CU(\cdot, \cdot)$ is defined in \eqref{def CU}.
Note that
 the third line is due to  the following identity:
\beqn
w_1'\circ f = (\CL-I)\wrho^{(1)}\circ f
=\CL\wrho^{(1)}\circ f - \wrho^{(1)}\circ f 
= \wrho^{(1)}_{\mathrm{even}}  - \wrho^{(1)}\circ f,
\eeqn
and the fourth line uses that 
$\ds
w_1=c_1+\ww_1
$.
To summarize, by Proposition~\ref{prop: qu req}, 
we have
\beqn
\begin{cases}
v_1= \xi_1 + w_1\circ f =c_1 + \xi_1 +   \hv_1, \\
v_2= \xi_2 + w_2\circ f + w_{2\#} \circ f  = c_2 + \xi_2 +  \hv_2 +
\CU(c_1, \xi_1).
\end{cases}
\eeqn
The proof of this proposition is complete.
\end{proof}

\begin{remark}
\label{rem: full C5 solutions}
The preceding propositions describe a convenient smoother family of
solutions. For the full $C^5$ solution set in
Proposition~\ref{prop: qu req}, one may instead take
$\xi_1,\xi_2\in\rho^{-1}\ker(\CL|_{C^5(\TT)})$ and require
$w_{2\#}\in C^5(\TT)$. This condition is necessary by
\eqref{solve v3}, and sufficient by \eqref{sol qu req}.
For the doubling map, the full solution set is therefore given by the
same displayed formulas with
$\xi_1,\xi_2\in C^5(\TT)\cap\CB_{\mathrm{odd}}$ and
$\CU(c_1,\xi_1)\in C^5(\TT)$, equivalently
$\xi_1\xi_1'\in C^5(\TT)$.
\end{remark}

\medskip

Under the setting of Proposition~\ref{prop: qu req double}, for any function $\CJ: C^5(\TT)\times C^5(\TT)\to\RR$ satisfying conditions in
Corollary~\ref{cor: request opt},
there must exist a minimal solution 
$\ds (\ov_1,\ov_2)\in\overline{\CV}^{\,\beta}$ ($3<\beta<5$),
where $\CV:=\CV\left(f;\wrho^{(1)},\wrho^{(2)}\right)$ is the solution set
to the quadratic request problem,
and the closure $\overline{\CV}^{\,\beta}$ and extension $\overline{\CJ}$
are as defined in Corollary~\ref{cor: request opt}, such that
\beqn
\overline{\CJ}(\ov_1,\ov_2)
=\inf\left\{\CJ(v_1,v_2):\ (v_1,v_2)\in \CV\left(f;\wrho^{(1)},\wrho^{(2)}\right)\right\}.
\eeqn
The following example shows that the minimal solutions need not be unique.

\begin{example}
\label{example: quadratic non-uniq}
Under the setting of Proposition~\ref{prop: qu req double}, i.e.,  the unperturbed map $f$ is 
the doubling map of $\TT$, we further consider the special prescribed responses $\ds \wrho^{(1)}(x)\equiv 0$ and  $\ds \wrho^{(2)}(x)= \pi^2\tau^2\cos(2\pi x) $ for some $\tau\in \RR$. 
Set $\alpha=5$, fix $3<\beta<5$, and write $\CV:=\CV^5\left(f;\wrho^{(1)},\wrho^{(2)}\right)$.
By Proposition~\ref{prop: qu req double}, Remark~\ref{rem: full C5 solutions},
and direct computation, we obtain the pairs in $\CV$ are precisely those of the form
\beqn
\begin{cases}
v_1= c_1 + \xi_1  , \\
v_2= c_2 + \xi_2 + \frac12\xi_1\xi_1' +  \hv_2,    
\end{cases}
\eeqn
where $c_1, c_2\in \RR$, $\ds\xi_1\in C^{5}(\TT)\cap\CB_{\mathrm{odd}}$,
$\ds\xi_2\in C^5(\TT)\cap\CB_{\mathrm{odd}}$,
with $\xi_1\xi_1'\in C^5(\TT)$, and $\ds \hv_2(x)= -\frac12 \pi \tau^2 \sin(4\pi x) \in \CB_{\mathrm{even}}$. 
Let us further  consider the particular function $\CJ_\lambda:C^5(\TT)\times C^5(\TT)\to\RR$ for $\lambda>0$, given by
\beqn
\CJ_\lambda(v_1,v_2)
:=\|v_1\|_{L^2}^2+\|v_2\|_{L^2}^2
+\lambda\left(\|v_1\|_{C^5}^2+\|v_2\|_{C^5}^2\right),
\eeqn
where $\|u\|_{C^5}:=\sum_{j=0}^5\|u^{(j)}\|_\infty$.
The function $\CJ_\lambda$ satisfies conditions (1) and (2) of
Corollary~\ref{cor: request opt}.
The $C^5$ penalty supplies the coercivity absent from the unregularized $L^2$ cost. Its lower semicontinuity under $C^\beta$ convergence to a $C^5$ limit follows from distributional convergence of derivatives and lower semicontinuity of their supremum norms.

Then Corollary~\ref{cor: request opt} implies that there exists a minimal solution $(\ov_1,\ov_2)\in\overline{\CV}^{\,\beta}$ with respect to the extension $\overline{\CJ_\lambda}$.
We claim that 
the minimal solution need  not  be unique. 
Indeed, passing to the limit in the solution formula, any $(0,v_2)\in\overline{\CV}^{\,\beta}$ has the form $v_2=c_2+\xi_2+\hv_2$, where $c_2\in\RR$ and $\xi_2\in\CB_{\mathrm{odd}}$.
By  $L^2$-orthogonality between $\CB_{\mathrm{odd}}$ and $\CB_{\mathrm{even}}$, together with $\int_\TT\hv_2(x)\,dx=0$, we have
\beqn
\begin{split}
\overline{\CJ_\lambda}(0,v_2)
&\ge\|v_2\|_{L^2}^2
=c_2^2+\|\xi_2\|_{L^2}^2+\|\hv_2\|_{L^2}^2\\
&\ge\|\hv_2\|_{L^2}^2=\pi^2\tau^4/8.
\end{split}
\eeqn 
On the other hand, let $\xi_1^*(x)=\tau\sin(2\pi x)$, then $(\xi_1^*,0)\in\CV$ as  $\frac12\xi_1^*(\xi_1^*)'=-\hv_2$. Thus,
	\beqn
	\overline{\CJ_\lambda}(\xi_1^*,0)
	=\CJ_\lambda(\xi_1^*,0)
	=\tau^2\left(\frac12+\lambda\|\sin(2\pi x)\|_{C^5}^2\right).
	\eeqn
Therefore, if $\ds\tau>2/\sqrt\pi$ and
	$0<\lambda<(\pi^2\tau^2-4)/(8\|\sin(2\pi x)\|_{C^5}^2)$, then
\beqn
\overline{\CJ_\lambda}(\xi_1^*,0)<\pi^2\tau^4/8.
\eeqn
Thus, if $(\ov_1,\ov_2)$ is a minimal solution for $\overline{\CJ_\lambda}$, we must have $\ov_1\not\equiv0$. Note $(-v_1,v_2)\in\overline{\CV}^{\,\beta}$ and $\overline{\CJ_\lambda}(-v_1,v_2)=\overline{\CJ_\lambda}(v_1,v_2)$ for any $(v_1,v_2)\in\overline{\CV}^{\,\beta}$.
Therefore, we get that $(-\ov_1,\ov_2)$  also achieves the minimum of $\overline{\CJ_\lambda}$.
Therefore,  if $\ds \tau> 2/\sqrt \pi$ and 
$0<\lambda<(\pi^2\tau^2-4)/(8\|\sin(2\pi x)\|_{C^5}^2)$, there are at least two distinct minimizers of
$\overline{\CJ_\lambda}$ on $\overline{\CV}^{\,\beta}$.
For the explicit family $f_t(x)=2x+t\xi_1^*(x)\pmod1$, the observable $\cos(2\pi x)$ has expectation
\[
\int_{\TT}\cos(2\pi x)\rho_t(x)\,dx
=\frac{\pi^2\tau^2}{2}t^2+o(t^2).
\]
Thus a perturbation invisible to linear response still changes a concrete statistic at second order. The symmetry $v_1\mapsto-v_1$ explains the two distinct relaxed optima.
\end{example}

\bigskip  
  
\section*{Acknowledgment}
The authors thank the anonymous referees for their careful reading and constructive comments, which helped improve the geometric formulation, the perturbation argument, and the presentation of this paper.
J. Chen is partially supported by the National Key Research and Development Program of China (No.2022YFA1005802), 
the NSFC Grant 12471186,
NSF of Jiangsu BK20241916
and the Jiangsu Shuang chuang Project (JSSCTD202209).

\bibliographystyle{plain}
\bibliography{reference1}{}

\end{document}